\documentclass{amsart}
\usepackage{a4wide}
\usepackage[latin1]{inputenc}
\usepackage[T1]{fontenc}
\usepackage{amsmath,amssymb,amsthm,stmaryrd}
\usepackage{mathrsfs}
\usepackage{esint}
\usepackage{color}
\usepackage{pdfsync}

\newcommand{\ud}[0]{\,\mathrm{d}}

\newcommand{\eps}[0]{\varepsilon}

\newcommand{\loc}[0]{\operatorname{loc}}
\newcommand{\diam}[0]{\operatorname{diam}}

\newcommand{\abs}[1]{|#1|}
\newcommand{\Babs}[1]{\Big|#1\Big|}

\newcommand{\Norm}[2]{\|#1\|_{#2}}
\newcommand{\BNorm}[2]{\Big\|#1\Big\|_{#2}}
\newcommand{\pair}[2]{\langle #1,#2 \rangle}

\newcommand{\lspan}[0]{\operatorname{span}}

\newcommand{\BMO}[0]{\operatorname{BMO}}
\newcommand{\supp}[0]{\operatorname{supp}}

\newcommand\R{\mathbb{R}}
\newcommand\C{\mathbb{C}}
\newcommand\N{\mathbb{N}}
\newcommand\Z{\mathbb{Z}}

\newcommand\Q{\mathbb{Q}}

\newcommand{\prob}[0]{\mathbb{P}}

\swapnumbers

\numberwithin{equation}{section}

\makeatletter
  \let\c@subsection\c@equation
\makeatother

\theoremstyle{plain}
\newtheorem{theorem}[equation]{Theorem}
\newtheorem{proposition}[equation]{Proposition}
\newtheorem{corollary}[equation]{Corollary}
\newtheorem{lemma}[equation]{Lemma}

\theoremstyle{definition}
\newtheorem{definition}[equation]{Definition}

\theoremstyle{remark}
\newtheorem{remark}[equation]{Remark}

\makeatletter
\@namedef{subjclassname@2010}{
  \textup{2010} Mathematics Subject Classification}
\makeatother

\title[Wavelets in spaces of homogeneous type]{Orthonormal bases of regular wavelets in spaces of homogeneous type}
\author{Pascal Auscher}
\address{Univ. Paris-Sud, laboratoire de Math\'ematiques, UMR 8628, F-91405 {\sc Orsay}; CNRS, F-91405 {\sc Orsay}} 
\email{pascal.auscher@math.u-psud.fr}
\author{Tuomas Hyt\"onen}
\address{Department of Mathematics and Statistics, P.O. Box 68, FI-00014 University of Helsinki, Finland}
\email{tuomas.hytonen@helsinki.fi}


\begin{document}

\begin{abstract} Adapting the recently developed randomized dyadic structures, we introduce the notion of spline function in geometrically doubling quasi-metric spaces. 
Such functions have interpolation and reproducing properties  as the linear splines in Euclidean spaces. They also have H\"older regularity. This is used to build an orthonormal basis of H\"older-continuous wavelets with exponential decay in any space of homogeneous type. As in the classical theory, wavelet bases provide a universal  Calder\'on reproducing formula  to study and develop function space theory and singular integrals. We discuss the examples of $L^p$ spaces, BMO and apply  this to a proof of the $T(1)$ theorem. As no extra condition {(like `reverse doubling', `small boundary' of balls, etc.)} on the space of homogeneous type is required,  our results extend a long line of works on the subject.  
\end{abstract}

\maketitle

\section{Introduction}

The main goals of this paper are two-fold: the construction of orthonormal, H\"older-continuous wavelet bases in general spaces of homogeneous type, and their applications in the theory of singular integrals and function spaces in the same general set-up. Despite several existing results of related nature, the scope of our theory is completely new in at least two respects: 

First, as far as we are aware, it seems that we offer the first construction of an orthonormal wavelet \emph{basis}, as opposed to a \emph{frame}, in this setting. Only relatively recently, Deng and Han wrote \cite[p.~40]{DH}: ``Orthonormal wavelet bases are out of reach on a space of homogeneous type. Instead the theory of frames will be used. Roughly speaking using a frame means that you tolerate a limited amount of redundancy while redundancy is completely avoided with a basis.'' Here, we build a genuine basis, and so avoid this redundancy. Our construction starts from appropriate \emph{splines}, which also seem to be new on an abstract quasi-metric, even metric,  space, and of independent interest.

Second, we are careful not to impose any additional assumptions other than those defining a space of homogeneous type in the sense of Coifman and Weiss \cite{CW}---that is, a set $X$ equipped with a quasi-distance $d$ and a  Borel measure $\mu$ that is doubling on the quasi-metric balls---, and we keep working with the original given $d$, without changing to an `equivalent' one (see below).  The difficulty is that  the quasi-distance, in contrast to a distance,  may not be H\"older-regular, and quasi-metric balls might not be open, nor even Borel sets with respect to the topology defined by the quasi-distance (openness is assumed in \cite{CW}).  We merely assume that  balls be Borel sets,  but  even this assumption can be easily remedied by. See  Section \ref{sec:l2}.

In contrast to this, we emphasize that a substantial part of the rich literature on analysis in 
spaces of homogeneous type
 involves additional assumptions on the underlying space. {Let us consider in particular the existence of a Littlewood--Paley decomposition (or a Calder\'on reproducing formula) with a \emph{regular} kernel, which has been the object of numerous works. (Decompositions with discontinuous kernels arising from martingale differences are well-known even in more general set-ups; this is a different story.) Indeed, the regular Littlewood--Paley decomposition is the core formula that allows the various characterizations of Hardy spaces, the development of function spaces  of Besov and Triebel--Lizorkin type, and the analysis of operators that act on them.
(Here, we have in mind a definition of Besov spaces linked with the modulus of continuity and thereby the (quasi-)metric of the underlying space; some other generalizations from the point-of-view of abstract approximation spaces, as in \cite{MM}, are yet another story.)

As far as we understand, there has always been some kind of further hypothesis either on $X$, $d$ or $\mu$. The seminal work in this context is the paper of David--Journ\'e--Semmes on the $T(b)$ theorem \cite{DJS}. The starting point of their theory is the existence of a new H\"older-continuous quasi-metric {adapted to the measure $\mu$ and}
\emph{topologically} (but not geometrically) equivalent to the original $d$, which is provided by work of Macias--Segovia \cite[{Theorems 2 and 3}]{MS} { (for more precise {versions} of the Macias-Segovia results, see  \cite{PS} and \cite{MMMM})}. Their singular integrals are then defined relative to this new quasi-metric, which leads to a {possibly} different class of Calder\'on--Zygmund operators than those defined in terms of the original $d$ and $\mu$. ({{Even then, \cite{DJS} still needs}
some more technical assumptions such as no point masses, infinite mass for $X$ and even a ``small boundary'' property for balls.)
Working with this kind of new quasi-metric has become a common set-up in the literature.
The theory of function spaces under this assumption has been developed by Han--Sawyer \cite{HS}, see also the nice review by Han--Weiss \cite{HW}.

{ By \cite{MS}, Theorem 2 (for an elementary proof, see \cite{PS}), it is possible to  only change the quasi-metric to a metrically equivalent one that is H\"older-regular. This does not change the classes of singular integral operators (up to changing the H\"older-exponents), but  leaves another difficulty:  some estimates require a control from below for the growth of balls.  See the discussion below about \eqref{eq:technical}. 
In several more recent papers \cite{HMY2,HMY1,YZ2,YZ1} on function spaces on (quasi-)metric spaces, this issue is circumvented by assuming}
 the  `reverse doubling' property
\begin{equation*}
  \mu(B(x,Cr))\geq(1+\varepsilon)\mu(B(x,r)),\qquad 0<r<\diam(X)/C,
\end{equation*}
which is equivalent to the non-emptyness of annuli: $B(x,R)\setminus B(x,r)\neq\varnothing$ whenever $R/r$ is large enough and $r$ is at most a fraction of $\diam(X)$. This latter condition is another common assumption;  it was only recently eliminated from some results about positive integral operators on metric spaces by Kairema \cite{Kairema}.

The reverse doubling excludes in particular the presence of point masses (which is also frequently assumed, even if empty annuli are otherwise tolerated), and therefore rules out some basic examples like discrete groups  $\Z$ or $\Z/p\Z$ and the multidimensional analogs, the typical discrete metric structures arising in theoretical computer science (trees, graphs, or strings from a finite alphabet), $\Q_{p}$ from arithmetics, or discrete approximations of other spaces of homogeneous type (like those constructed in \cite{ACI}), even if the original space did satisfy the reverse doubling property and some \textit{ad hoc} orthonormal (or bi-orthogonal) wavelet bases had been constructed. To mention a few references in these directions, see \cite{AES}, \cite{PZ}, \cite{S}. We remark that even if H\"older-regularity on discrete structures might not be an issue at first glance, it becomes one thinking the discrete structure as approximating a continuous one.

In contrast, we repeat, we do not make any such additional assumptions or changes to a different quasi-metric.

Let us come to the construction of spline wavelets.  Splines have a long history in approximation theory, due in particular to the interpolation property  and   their polynomial nature, and were popularized by early books in the 1960's. (The first book listed from a MathSciNet title search ``spline'' is \cite{ANW} and the reviewer dates the appearance of the splines  to Schoenberg's work  during World War II; see \cite{Sch}.)  As for spline wavelets, they were first constructed on the real line and on the 1-torus by Str\"omberg \cite{Str}. They were rediscovered  independently by Battle \cite{Bat} and Lemari\'e \cite{Lem1} in Euclidean spaces. Using the Multiresolution Analysis scheme of Mallat \cite{Mal} and Meyer \cite{Mey}, it became clear that the existence of regular, compactly supported splines with the interpolation property and reproducing formula  leads to bases of regular spline wavelets. This was attempted in other geometrical contexts. For example, splines  and spline wavelets were constructed in \cite{Lem2} on a stratified Lie group following a minimisation procedure using the existence of a sublaplacian, and in \cite{DDW} on compact manifolds like spheres. We do not claim to be exhaustive in this history. Nevertheless,  we have to  follow a new route as  there is no group structure attached, nor any local coordinates to help. We rather use a probabilistic approach as described next.

Central to our analysis are the dyadic structures in a geometrically doubling quasi-metric space  as  constructed by Christ \cite{Christ} (a construction in Ahlfors-David spaces was done earlier by David \cite{D}),  and their randomized versions recently studied by the second author with  Martikainen \cite{HM} and Kairema \cite{HK}. {(Another variant is due to Nazarov, Reznikov and Volberg \cite{NRV}.)} We build the splines as averages of the indicators of dyadic cubes under a random choice of the dyadic system. This can be motivated by an appropriate point of view at  the classical piecewise linear splines on $\R$, which are generated by the function
\begin{equation*}
  s(x)= x1_{(0,1]}(x)+(2-x)1_{(1,2)}(x).
\end{equation*}
Observe that random dyadic intervals of sidelength $1$, in the sense of Nazarov, Treil and Volberg {\cite[Sec.~9.1]{NTV}}, can be defined by translating the standard intervals $[k,k+1)$ by a random number $u\in[0,1)$. Thus the unit cube with left end at the origin, $[0,1)$, is translated to $[u,u+1)$. It is easy to check that the average of its indicator function $1_{[u,u+1)}(x)=1_{[0,1)}(x-u)$ over the uniform random choice of $u$ is precisely
\begin{equation*}
  \int_0^1 1_{[0,1)}(x-u)\ud u =1_{[0,1)}*1_{[0,1)}(x) = s(x).
\end{equation*}
So one can indeed think of splines as averages of the indicators of random dyadic intervals. This is the basic idea which will guide us in constructing splines for a quasi-metric space $X$.

Once the splines are available, and only at this point we add into our considerations the doubling Borel measure, we can obtain the orthonormal bases of scaling functions and the wavelets by the general procedure of Meyer. These will be H\"older-continuous functions $\psi^k_\alpha$ with an  \emph{exponential localization} in the $\delta^k$-neighbourhood of a given point $y^k_\alpha$:
\begin{equation*}
  \abs{\psi^k_\alpha(x)}\leq\frac{C}{\sqrt{\mu(B(y^k_\alpha,\delta^k))}}\exp\Big(-\gamma \Big(\frac{d(x,y^k_\alpha)}{\delta ^{k}}\Big)^a\, \Big)
\end{equation*}
for some $a\in (0,1]$ that depends only on the  quasi-triangle constant (with $a=1$ in the metric case or when the quasi-distance is Lipschitz-continuous).
The construction of \emph{boundedly supported} wavelets in this generality remains an open problem.

There is an important feature related to the collections of points $\mathscr{Y}^k:=\{y^k_\alpha\}$, which index the wavelets of the scale $\delta^k$. These collections arise as $\mathscr{Y}^k=\mathscr{X}^{k+1}\setminus\mathscr{X}^k$, where $\cdots\subseteq\mathscr{X}^k\subseteq\mathscr{X}^{k+1}\subseteq\cdots$ is an increasing sequence of point-sets, and each $\mathscr{X}^k$ is both $\delta^k$-separated and (up to constant) $\delta^k$-dense in $X$. The important feature is that $\mathscr{Y}^k$ can be much sparser than $\mathscr{X}^{k+1}$ in some regions of the space, which reflects the `absence of the length scale $\delta^{k+1}$' in the local geometry of the space. Accordingly, the distance $d(x,\mathscr{Y}^k)$ of a given element $x\in X$ to the nearest point $y^k_\alpha\in\mathscr{Y}^k$ will be a significant quantity, which can in general be much larger that $\delta^k$.

This quantity appears in central technical estimates when bounding series of the following type, which naturally arise in the context of function spaces:
\begin{equation}\label{eq:technical}
  \sum_{j\in\Z:\delta^j\geq r}\frac{1}{\mu(B(x,\delta^j))}\exp\Big(-\gamma\Big(\frac{d(x,\mathscr{Y}^j)}{\delta^j}\Big)^a\, \Big)\leq\frac{C}{\mu(B(x,r))}.
\end{equation}
This (non-obvious but valid) estimate serves as a replacement of the following bound, which is repeatedly applied by Han, M\"uller and Yang in the reverse doubling context (cf. \cite[Lemma 3.5]{HMY1}):
\begin{equation*}
  \sum_{j\in\Z:\delta^j\geq r}\frac{1}{\mu(B(x,\delta^j))}\leq\frac{C}{\mu(B(x,r))}.
\end{equation*}
Indeed, this is a quick consequence of reverse doubling, but invalid in general spaces of homogeneous type. A typical cause to destroy this latter estimate is the presence of a large empty annulus around $x$, but this then results in some large values of $d(x,\mathscr{Y}^j)$, which compensates for the failure of $\mu(B(x,\delta^j))$ to grow fast enough in \eqref{eq:technical}.

Although we do not develop the theory of function spaces in the detail of \cite{HMY1} here, it is clear that our analysis involving the quantity $d(x,\mathscr{Y}^k)$ will extend to further related questions that we do not explicitly deal with; we believe that it could be used to extend large parts of the recent theory of `RD (reverse doubling) spaces' to general spaces of homogeneous type.
As an illustration,  we give a short  treatment the $T(1)$ theorem using our  construction.

{Over the last decade or so, it has been discovered that several aspects of harmonic analysis can even be pushed beyond the setting of doubling measures; see e.g. \cite{NTV,Tolsa:BMO,Verdera} for some of the pioneering developments in the context of $\R^d$, and \cite{HLYY,HM} for recent developments in abstract quasi-metric spaces.}
We do not address the non-doubling measures here, as it seems that our regular splines and wavelets are most useful in the (already quite general) doubling situation. In fact, while our splines are constructed completely independently of any underlying measure, it turns out that they automatically have a good $L^2$ theory with respect to any doubling measure on the space. This would not be the case for a non-doubling measure, and a successful wavelet theory for such measures, if possible at all, should probably be adapted to the particular measure in a more complicated manner, already at the early point, where we can manage with a purely geometric construction.

\subsection*{About notation}
We use $C$ to denote positive constants, whose value may change from one occurrence to the next. We also abbreviate `$\leq C\times\cdots$' to just `$\lesssim$'. We use $\gamma$ in a similar role as a positive exponent. In contrast to $C$, it typically decreases from one occurrence to the next. For the measure (`volume') of balls {$B(x,r):=\{y\in X;d(x,y)<r\}$}, we sometimes use the abbreviations
\begin{equation*}
  V(x,r):=\mu(B(x,r)),\qquad V(x,y):=V(x,d(x,y)).
\end{equation*}

\subsection*{Acknowledgements} This work started during a visit of the second author to Paris-Sud whose support is gratefully acknowledged. {The second author was also partially supported by the Academy of Finland, grants 130166, 133264 and 218148.} We thank Detlef M\"uller for answering our questions on his work on function spaces on Reverse Doubling spaces, and the anonymous referee for comments to improve the presentation.

\section{Preliminaries on dyadic systems}

\subsection*{Introduction and motivation}
In this section we review the relevant-for-us parts of M.~Christ's \cite{Christ} construction of dyadic cubes in a space of homogeneous type, as well as its recent probabilistic version from \cite{HK,HM}. For the convenience of the reader, and since we feel that we have managed to slightly streamline the earlier presentations, we will give a self-contained treatment, even though the actual novelty is only in certain details of the construction.

Christ's cubes $Q^k_\alpha$ (where $k$ designates the generation or length scale of the cube, and $\alpha$ indexes the cubes inside a given generation) are defined by two sets of auxiliary objects: the centre points $x^k_\alpha$, which determine the rough position of individual cubes, and a partial order $\leq$ on the family of the index pairs $(k,\alpha)$, which determines the set inclusions among different cubes.

In \cite{HK,HM}, a systematic method of constructing several Christ-type families of dyadic cubes $Q^k_\alpha(\omega)$ was introduced, where $\omega$ belongs to a parameter space $\Omega$. Here $\Omega$ can be equipped with a probability measure, which gives rise to a notion of random dyadic cubes. This allows us to compute averages over a random choice of the cubes, and make probabilistic statements about them. As usual, averaging has a smoothing effect, and the random cubes will `on average' enjoy better regularity properties than any fixed cubes would do. This is essential for our construction of H\"older-continuous splines.

Since the cubes $Q^k_\alpha$ are determined by the dyadic points $x^k_\alpha$ and the partial order $\leq$, the construction of the parametrized family of $Q^k_\alpha(\omega)$ amounts to defining appropriate parameter-dependent points $z^k_\alpha(\omega)$ and a parameter-dependent partial order $\leq_\omega$. This is achieved by first fixing the reference objects $x^k_\alpha$ and $\leq$ as in Christ's original work, and constructing $z^k_\alpha(\omega)$ and $\leq_\omega$ as their parametrized perturbations.

All these constructions involve certain arbitrary choices, which all work equally well for the applications considered in the earlier papers. A feature of our presentation here is that we are going to insist on somewhat more specific, rather than arbitrary, choices in certain details of the construction, as this will be a convenience later when using the dyadic cubes in the building of our splines.

We now turn to the details.

\subsection*{General assumptions}
In what follows, $X$ is a set equipped with a quasi-distance with quasi-triangle constant $A_{0}\geq 1$, namely, $d(x,y)=d(y,x)\ge 0$ , $d(x,y)=0$ if and only if $x=y$, and {$d$ satisfies the} quasi-triangle  inequality
$$
d(x,y) \le A_{0} (d(x,z)+d(z,y)).
$$
We assume that $X$ has \emph{the geometric doubling property}, namely that there exists a natural number $N$ such that  any given ball  contains no more than $N$ points at quasi-distance exceeding half its radius.
This is a certain finite-dimensionality requirement on the space: for example, the infinite-dimensional Hilbert space $\ell^2$ fails this property, since all the half-unit vectors $\tfrac12 e_k$, $k\in\N$, belong to the unit-ball, while their mutual distance is $\sqrt{2}/2>1/2$.
We do not need the doubling measure, or in fact any underlying measure, at this stage.

Let us fix a small parameter $\delta>0$. For example, it suffices to take $\delta\leq\frac{1}{1000}A_0^{-10}$.
Roughly speaking, the point is that phenomena on the length scale $\delta^{k+1}$ should remain much smaller than the length scale $\delta^k$, even after repeated use of the quasi-triangle inequality where we `lose' the constant $A_0$ at every application. For example, if the points $x_0,x_1,x_2\ldots$ satisfy $d(x_{i-1},x_i)<\delta^{k+1}$, we can still conclude that $d(x_0,x_r)$ is much less that $\delta^k$ for $r\leq 10$.

\subsection*{Reference dyadic points}
For every $k\in\Z$, we choose a set of \emph{reference dyadic points} $x^k_{\alpha}$ as follows: For $k=0$, let $\mathscr{X}^0:=\{x^0_{\alpha}\}_{\alpha}$ be a maximal collection of $1$-separated points. Inductively, for $k\in\Z_+$, let $\mathscr{X}^k:=\{x^{k}_{\alpha}\}_{\alpha}\supseteq\mathscr{X}^{k-1}$ and $\mathscr{X}^{-k}:=\{x^{-k}_{\alpha}\}_{\alpha}\subseteq\mathscr{X}^{-(k-1)}$ be maximal $\delta^{k}$- and $\delta^{-k}$-separated collections in $X$ and in $\mathscr{X}^{-(k-1)}$, respectively.

\begin{lemma}
For all $k\in\Z$ and $x\in X$, the reference dyadic points satisfy
\begin{equation*}
  d(x^k_{\alpha},x^k_{\beta})\geq\delta^k\quad\big(\alpha\neq\beta\big),\qquad d(x,\mathscr{X}^k)=\min_{\alpha}d(x,x^k_{\alpha})<2A_0\delta^k.
\end{equation*}
\end{lemma}

\begin{proof}
The separation property is part of the construction.
By maximality, it follows that for all $x\in X$ and $k\geq 0$,
\begin{equation*}
  d(x,\mathscr{X}^k)=\min_{\alpha}d(x,x^k_{\alpha})<\delta^k.
\end{equation*}
Also, given $x\in X$, we can recursively find points $x^0_{\alpha_0},x^{-1}_{\alpha_1},\ldots,x^{-k}_{\alpha_k}$ such that
\begin{equation*}
  d(x,x^0_{\alpha_0})<1,\quad d(x^0_{\alpha_0},x^{-1}_{\alpha_1})<\delta^{-1},\quad\ldots\quad d(x^{-(k-1)}_{\alpha_{k-1}},x^{-k}_{\alpha_k})<\delta^{-k},
\end{equation*}
and hence
\begin{equation*}
  d(x,x^{-k}_{\alpha_k})<\sum_{j=0}^k A_0^{j+1}\delta^{-(k-j)}<\frac{A_0\delta^{-k}}{1-A_0\delta}\leq 2A_0\delta^{-k}.\qedhere
\end{equation*}
\end{proof}

Note that $\mathscr{X}^k\subseteq\mathscr{X}^{k+1}$, so that every $x^k_{\alpha}$ is also a point of the form $x^{k+1}_{\beta}$, and thus of all the finer levels.
We denote $\mathscr{Y}^k:=\mathscr{X}^{k+1}\setminus\mathscr{X}^k$, and relabel these points as $\mathscr{Y}^k=\{y^k_\alpha\}_{\alpha}$. These points will play an important role as a parameter set of our wavelets, to be constructed.

\subsection*{Reference partial order}
We set up a \emph{partial order} $\leq$ among the pairs $(k,\alpha)$ as follows: Each $(k+1,\beta)$ satisfies $(k+1,\beta)\leq(k,\alpha)$ for exactly one $(k,\alpha)$, in such a way that
\begin{equation}\label{eq:leqProperties}
 d(x^{k+1}_{\beta},x^k_{\alpha})<\frac{1}{2A_0}\delta^k\quad\Longrightarrow\quad(k+1,\beta)\leq(k,\alpha)
 \quad\Longrightarrow\quad d(x^{k+1}_{\beta},x^k_{\alpha})<2A_0\delta^k.
\end{equation}
The pairs $(k+1,\beta)$ with $(k+1,\beta)\leq(k,\alpha)$ are called the children of $(k,\alpha)$. Geometric doubling implies that their number is uniformly bounded. So far, we have essentially followed the original construction of M.~Christ \cite{Christ}, a slight nuance being the choice of the point-sets $\mathscr{X}^k$ in such a way as to have the nestedness $\mathscr{X}^k\subseteq\mathscr{X}^{k+1}$, which was not required in \cite{Christ,HK,HM}.

\subsection*{Labels for the points}
For a successful perturbation argument to construct the parametrized dyadic points $z^k_\alpha(\omega)$ below, we need certain book-keeping among the near-by dyadic points $x^k_\alpha$ of the same generation.

Points $(k,\alpha)$ and $(k,\beta)$ are called \emph{neighbours}, if they have children $(k+1,\gamma)\leq(k,\alpha)$ and $(k+1,\eta)\leq(k,\beta)$ such that $d(x^{k+1}_{\gamma},x^{k+1}_{\eta})<(2A_0)^{-1}\delta^k$. In this case
\begin{equation*}
\begin{split}
 d(x^k_{\alpha},x^k_{\beta})
 &<A_0 d(x^k_{\alpha},x^{k+1}_{\gamma})
 +A_0^2 d(x^{k+1}_{\gamma},x^{k+1}_{\eta})
 +A_0^2 d(x^{k+1}_{\eta},x^k_{\beta}) \\
 &<2A_0^2\delta^k+\tfrac12 A_0\delta^k+2A_0^3\delta^k<5A_0^3\delta^k.
\end{split}
\end{equation*}
The number of neighbours that any point can have is also uniformly bounded.

We equip each pair $(k,\alpha)$ with two \emph{labels}, which are chosen from a finite set but which still locally distinguish between different pairs $(k,\alpha)$. The primary label, denoted by $\operatorname{label}_1(k,\alpha)\in\{0,1,\ldots,L\}$, where $L$ is the maximal number of neighbours, is chosen in such a way that any two neighbours have a different label. The secondary label, denoted by $\operatorname{label}_2(k,\alpha)\in\{1,\ldots,M\}$, where $M$ is the maximal number of children, is chosen in such a way that no two children of the same parent have the same label. These labels were introduced in \cite{HK} with slightly different notation.

\subsection*{Parametrized points and partial order}
As described above, we now want to perform a perturbation of the original $x^k_\alpha$ and $\leq$ so as to obtain a parametrized family of similar objects, on which probabilistic statements can later be made. The parameter space will be
\begin{equation*}
  \Omega=\Big(\{0,1,\ldots,L\}\times\{1,\ldots,M\}\Big)^{\Z},
\end{equation*}
with a typical point denoted by $\omega=(\omega_k)_{k\in\Z}$, where $\omega_k=(\ell_k,m_k)\in\{0,1,\ldots,L\}\times\{1,\ldots,M\}$.

The \emph{new dyadic points} $z^k_{\alpha}=z^k_{\alpha}(\omega_k)$ are defined by
\begin{equation*}
  z^k_{\alpha}:=\begin{cases}  x^{k+1}_{\beta} & \text{if }\operatorname{label}_1(k,\alpha)=\ell_k,\text{ and } (k+1,\beta)\leq(k,\alpha),\text{ and }\operatorname{label}_2(k+1,\beta)=m_k, \\
    x^k_{\alpha} &  \text{if }\operatorname{label}_1(k,\alpha)\neq\ell_k,\text{ or } \not\exists(k+1,\beta)\leq(k,\alpha)
    \text{ such that }\operatorname{label}_2(k+1,\beta)=m_k. \\
    \end{cases}
\end{equation*}

A key feature of this definition is the following probabilistic statement, when $\Omega$ is equipped with the natural probability measure $\mathbb{P}_\omega$, which makes all coordinates $\omega_k=(\ell_k,m_k)$ are independent of each other and uniformly distributed over the finite set $\{0,1,\ldots,L\}\times\{1,\ldots,M\}$ (in other words, $\ell_k$ is uniformly distributed over $\{0,1,\ldots,L\}$, and $m_k$ over $\{1,\ldots,M\}$, independently).

\begin{lemma}\label{lem:prob}
Let $(k+1,\beta)\leq(k,\alpha)$ be fixed. Then
\begin{equation*}
  \mathbb{P}_\omega( z^k_\alpha=x^{k+1}_\beta )\geq\frac{1}{(L+1)M}.
\end{equation*}
\end{lemma}

In other words, every old point on the level $k+1$ has a positive (and bounded from below) probability of being a new point on the level $k$.

\begin{proof}
A sufficient condition that $ z^k_\alpha=x^{k+1}_\beta$ is that $\ell_k=\operatorname{label}_1(k,\alpha)$ (which has probability $1/(L+1)$) and that $m_k=\operatorname{label}_2(k+1,\beta)$ (which has probability $1/M$). Multiplying the probabilities of these independent events gives the claim.
\end{proof}

The point of the next result is that the new points behave qualitatively like the reference points, only with slightly weaker constants.

\begin{lemma}
The new points satisfy
\begin{equation*}
  d(z^k_{\alpha},z^k_{\beta})\geq\frac{1}{2A_0}\delta^k,\qquad\min_{\alpha}d(x,z^k_{\alpha})<4A_0^2\delta^k.
\end{equation*}
\end{lemma}

\begin{proof}
Consider the second bound first. Note that, in either of the two possibilities for the new point, we have $z^k_{\alpha}=x^{k+1}_{\beta}$ for some $(k+1,\beta)\leq(k,\alpha)$; in particular, $d(x^k_{\alpha},z^k_{\alpha})<2A_0\delta^k$, and thus
\begin{equation*}
  \min_{\alpha}d(x,z^k_{\alpha})<\min_{\alpha}[A_0 d(x,x^k_{\alpha})+A_0 d(x^k_{\alpha},z^k_{\alpha})]<4A_0^2\delta^k.
\end{equation*}

Let us then estimate $d(z^k_{\alpha},z^k_{\beta})$ from below for $\alpha\neq\beta$. If $(k,\alpha)$ and $(k,\beta)$ are not neighbours, then by definition this distance is at least $(2A_0)^{-1}\delta^k$. So suppose that these pairs are neighbours. Then they have different primary labels, and hence at least one of the new dyadic points, say $z^k_{\alpha}$, must satisfy $z^k_{\alpha}=x^k_{\alpha}$. On the other hand, if $(k+1,\eta)\leq(k,\beta)$, then $(k+1,\eta)\not\leq(k,\alpha)$, and we thus know that $d(x^{k+1}_{\eta},x^k_{\alpha})\geq(2A_0)^{-1}\delta^k$. But $z^k_{\beta}$ will be one of these points $x^{k+1}_{\eta}$; thus $d(z^k_{\beta},z^k_{\alpha})\geq(2A_0)^{-1}\delta^k$, as we claimed.
\end{proof}

The \emph{new partial order} $\leq_{\omega}$, $\omega=(\omega_k)_{k\in\Z}$, is set up as follows. We declare that
\begin{equation}\label{eq:newLeq}
   (k+1,\beta)\leq_{\omega}(k,\alpha)\quad\overset{\operatorname{def}}{\Longleftrightarrow}\quad
   \begin{cases} d(x^{k+1}_{\beta},z^k_{\alpha})<\tfrac14 A_0^{-2}\delta^k,\qquad \text{or} \\
     (k+1,\beta)\leq(k,\alpha) \text{ and } \not\exists\gamma: d(x^{k+1}_{\beta},z^k_{\gamma})<\tfrac14 A_0^{-2}\delta^k.
   \end{cases}
\end{equation}
In other words, to find the new parent of $(k+1,\beta)$ for the new partial order $\leq_{\omega}$, we first check whether the reference point $x^{k+1}_{\beta}$ is close (within distance $\tfrac14 A_0^{-2}\delta^k$) to some new dyadic point $z^k_{\alpha}$. If yes, then the corresponding $(k,\alpha)$ will be the new parent of $(k+1,\beta)$. If no such close point exists, then we simply use the original partial order $\leq$ to decide the parent of $(k+1,\beta)$.

\subsection*{Properties of the new points and order}

\begin{lemma}
For any given $k,\alpha,\beta$, the truth or falsity of the relation $(k+1,\beta)\leq_{\omega}(k,\alpha)$ depends only on the component $\omega_k$ of $\omega$.
\end{lemma}

\begin{proof}
In the defining conditions on the right of \eqref{eq:newLeq}, the only dependence on $\omega$ is via the new dyadic points $z^k_\alpha$, $z^k_\gamma$, and they depend only on $\omega_k$ by definition.
\end{proof}

This explicit definition of $\leq_\omega$ in terms of  the geometric configuration of the points and the original partial order $\leq$ is a novelty of our construction, where \cite{HK,HM} required a condition similar to \eqref{eq:leqProperties}, which only specifies the relation up to certain degrees of freedom. For us, a condition analogous to \eqref{eq:leqProperties} is a \emph{consequence} of the definition:

\begin{lemma}\label{lem:newLeqVsDist}
\begin{equation*}
 d(z^{k+1}_{\beta},z^k_{\alpha})<\tfrac15 A_0^{-3}\delta^k\quad\Longrightarrow\quad(k+1,\beta)\leq_{\omega}(k,\alpha)
 \quad\Longrightarrow\quad d(z^{k+1}_{\beta},z^k_{\alpha})<5A_0^3\delta^k.
\end{equation*}
\end{lemma}

\begin{proof}
Let first $d(z^{k+1}_{\beta},z^k_{\alpha})<\tfrac15 A_0^{-3}\delta^k$. Then
\begin{equation*}
  d(x^{k+1}_{\beta},z^k_{\alpha})
  \leq A_0d(x^{k+1}_{\beta},z^{k+1}_{\beta})+A_0d(z^{k+1}_{\beta},z^k_{\alpha})
  <2A_0^2\delta^{k+1}+\tfrac15 A_0^{-2}\delta^k\leq\tfrac14A_0^{-2}\delta^k,
\end{equation*}
and hence $(k+1,\beta)\leq_{\omega}(k,\alpha)$ by definition.

Let then $(k+1,\beta)\leq_{\omega}(k,\alpha)$. If $d(x^{k+1}_{\beta},z^k_{\alpha})<\tfrac14 A_0^{-2}\delta^k$, then
\begin{equation*}
  d(z^{k+1}_{\beta},z^k_{\alpha})\leq A_0d(z^{k+1}_{\beta},x^{k+1}_{\beta})+A_0d(x^{k+1}_{\beta},z^k_{\alpha})
  <2A_0^2\delta^{k+1}+\tfrac14A_0^{-1}\delta^k<5A_0^3\delta^k.
\end{equation*}
Otherwise, we have $(k+1,\beta)\leq(k,\alpha)$, and hence
\begin{equation*}
\begin{split}
  d(z^{k+1}_{\beta},z^k_{\alpha})
  &\leq A_0^2d(z^{k+1}_{\beta},x^{k+1}_{\beta})+A_0^2d(x^{k+1}_{\beta},x^k_{\alpha})+A_0d(x^k_{\alpha},z^k_{\alpha}) \\
  &<2A_0^3\delta^{k+1}+2A_0^3\delta^k+2A_0^2\delta^k< 5A_0^3\delta^k.\qedhere
\end{split}
\end{equation*}
\end{proof}

We can iterate this as follows:

\begin{lemma}\label{lem:newLeqVsDistIter}
For all $\ell\geq k$,
\begin{equation*}
 d(z^{\ell}_{\beta},z^k_{\alpha})<\tfrac16 A_0^{-4}\delta^k\quad\Longrightarrow\quad(\ell,\beta)\leq_{\omega}(k,\alpha)
 \quad\Longrightarrow\quad d(z^{\ell}_{\beta},z^k_{\alpha})<6A_0^4\delta^k.
\end{equation*}
\end{lemma}

\begin{proof}
The second implication follows from the second implication of Lemma~\ref{lem:newLeqVsDist} with the triangle inequality:
\begin{equation*}\label{eq:distOfDescendants}
 (\ell,\beta)\leq_{\omega}(k,\alpha)\quad\Longrightarrow\quad
  d(z^{\ell}_{\beta},z^k_{\alpha})<
  \sum_{j=0}^{\ell-k-1} 5A_0^3\delta^{k+j}\cdot A_0^{j+1}
  <\frac{5A_0^4\delta^k}{1-A_0\delta}<6A_0^4\delta^k.
\end{equation*}
For the first implication, if $\ell=k$, then the closeness of the points implies that $\beta=\alpha$. If $\ell>k$, consider $\gamma$ such that $(\ell,\beta)\leq_\omega(k+1,\gamma)$. By what we just proved, $d(z^\ell_\beta,z^{k+1}_\gamma)<6A_0^4\delta^{k+1}$, and hence
\begin{equation*}
  d(z^{k+1}_\gamma,z^k_\alpha)
  \leq A_0\Big(\frac{1}{6}A_0^{-4}\delta^k+6A_0^4\delta^{k+1}\Big)<\frac{1}{5}A_0^{-3}\delta^k;
\end{equation*}
thus $(\ell,\beta)\leq_\omega(k+1,\gamma)\leq_\omega(k,\alpha)$, where the last step follows from Lemma~\ref{lem:newLeqVsDist}.
\end{proof}

\subsection*{Dyadic cubes}
With the auxiliary objects at hand, the \textbf{dyadic cubes} are easy to define.
As in \cite{HK},  we introduce three families of these cubes---the preliminary, the closed, and the open:
\begin{equation*}
  \hat{Q}^k_{\alpha}(\omega):=\{z^{\ell}_{\beta}(\omega):(\ell,\beta)\leq_{\omega}(k,\alpha)\},
\end{equation*}
\begin{equation*}
  \bar{Q}^k_{\alpha}(\omega):=\overline{\hat{Q}^k_{\alpha}(\omega)},\qquad
  \tilde{Q}^k_{\alpha}(\omega):=\operatorname{interior}\bar{Q}^k_{\alpha}(\omega).
\end{equation*}
Note that $\hat{Q}^k_{\alpha}(\omega)$, and hence $\bar{Q}^k_{\alpha}(\omega)$ and $\tilde{Q}^k_{\alpha}(\omega)$, only depends on $\omega_{\ell}$ for $\ell\geq k$.

The rest of the section is concerned with the properties of these cubes. We first deal with several properties valid for an arbitrary fixed choice of the parameter $\omega\in\Omega$, and finally present a probabilistic statement concerning a random choice of $\omega$ under the natural product probability measure on the set $\Omega$.

\begin{lemma}\label{lem:cubeVsCentre}
\begin{equation*}
  \bar{Q}^k_{\alpha}(\omega)\subseteq B(z^k_{\alpha},6A_0^4\delta^k).
\end{equation*}
\end{lemma}

\begin{proof}
Let $x\in\bar{Q}^k_{\alpha}(\omega)$. Then $x$ is a limit of some points $z^{\ell}_{\beta}$ with $(\ell,\beta)\leq_{\omega}(k,\alpha)$, and we may choose a subsequence with $(\ell,\beta)\leq_{\omega}(k+1,\gamma)\leq_{\omega}(k,\alpha)$ for some $\gamma$. Then
\begin{equation*}
  d(z^k_{\alpha},x)
  \leq A_0 d(z^k_{\alpha},z^{k+1}_{\gamma})+A_0^2 d(z^{k+1}_{\gamma},z^{\ell}_{\beta})+ A_0^2 d(z^{\ell}_{\beta},x),
\end{equation*}
where the sum of the first two terms, by Lemma~\ref{lem:newLeqVsDist} and \eqref{eq:distOfDescendants} is $<5A_0^4 \delta^k+6A_0^6\delta^{k+1}<6A_0^4\delta^k$, and the last term becomes arbitrarily small.
\end{proof}

\begin{lemma}
We have the following covering properties for each fixed $k\in\Z$:
\begin{equation*}
  X=\bigcup_{\alpha}\bar{Q}^k_{\alpha}(\omega),\qquad
  \bar{Q}^k_{\alpha}(\omega)=\bigcup_{\beta:(k+1,\beta)\leq_{\omega}(k,\alpha)}\bar{Q}^{k+1}_{\beta}(\omega).
\end{equation*}
\end{lemma}

\begin{proof}
Every $x\in X$ is the limit of points $z^m_{\theta(m)}$ where $m\to\infty$. When $m\geq k$, we have $(m,\theta(m))\leq(k,\alpha)$ for some $\alpha$, and hence $\bigcup_\alpha\hat{Q}^k_\alpha(\omega)$ is dense in $X$. On the other hand, the uniform separation of the centres $z^k_\alpha$, the uniformly bounded diameter of the sets $\hat{Q}^k_\alpha(\omega)$, and the geometric doubling property imply that this union is locally finite. Hence the closure of the union is the union of the closures.

Similarly, it is immediate that $\hat{Q}^k_\alpha(\omega)=\bigcup_{\beta:(k+1,\beta)\leq_\omega(k,\alpha)}\hat{Q}^{k+1}_\beta(\omega)$, and the second claim follows by taking the closures, observing that the union is finite.
\end{proof}

\begin{lemma}
The closed and open cubes of the same generation are disjoint:
\begin{equation*}
  \bar{Q}^k_{\alpha}(\omega)\cap\tilde{Q}^k_{\beta}(\omega)=\varnothing\qquad(\alpha\neq\beta).
\end{equation*}
In fact,
\begin{equation*}
  \tilde{Q}^k_\alpha(\omega)^c=\bigcup_{\beta\neq\alpha}\bar{Q}^k_\beta(\omega).
\end{equation*}
\end{lemma}

\begin{proof}
We first check the weaker statement that
\begin{equation*}
   \bar{Q}^k_{\alpha}(\omega)\cap\hat{Q}^k_{\beta}(\omega)=\varnothing\qquad(\alpha\neq\beta).
\end{equation*}
Indeed, for contradiction, let $x\in\bar{Q}^k_{\alpha}(\omega)\cap\hat{Q}^k_{\beta}(\omega)$. Thus $x=z^\ell_\gamma$ with $(\ell,\gamma)\leq_\omega(k,\beta)$, and also $x=\lim_{m\to\infty}z^m_{\theta(m)}$ with $(m,\theta(m))\leq(k,\alpha)$. For $m$ large enough, we deduce that $d(z^m_{\theta(m)},z^\ell_\gamma)<\tfrac16A_0^{-4}\delta^\ell$ and thus, by Lemma~\ref{lem:newLeqVsDistIter}, $(m,\theta(m))\leq_\omega(\ell,\gamma)\leq_\omega(k,\beta)$, a contradiction.

To prove the actual first claim, again by contradiction, let $x\in \bar{Q}^k_{\alpha}(\omega)\cap\tilde{Q}^k_{\beta}(\omega)$. Hence $x=\lim_{m\to\infty}x_m$ with $x_m\in\hat{Q}^k_\alpha(\omega)$. Since $\tilde{Q}^k_{\beta}(\omega)$ is open by definition and $x_m\to x$, it follows that $x_m\in\tilde{Q}^k_\beta(\omega)\subseteq\bar{Q}^k_\beta(\omega)$ for large enough $m$. So in fact $x_m\in\hat{Q}^k_\alpha(\omega)\cap\bar{Q}^k_\beta(\omega)$, but this is empty by the first part of the proof. The claim follows.

For the second claim, it is immediate from the first claim that $\tilde{Q}^k_\alpha(\omega)\subseteq\bigcap_{\beta\neq\alpha}\bar{Q}^k_\beta(\omega)^c=\Big(\bigcup_{\beta\neq\alpha}\bar{Q}^k_\beta(\omega)\Big)^c=:\mathscr{O}^k_\alpha$. On the other hand, by local finiteness, the union defining $(\mathscr{O}^k_\alpha)^c$ is closed, and therefore $\mathscr{O}^k_\alpha$ itself, open. Since $X=\bar{Q}^k_\alpha(\omega)\cup(\mathscr{O}^k_\alpha)^c$, we have $\mathscr{O}^k_\alpha\subseteq\bar{Q}^k_\alpha(\omega)$, and since $\tilde{Q}^k_\alpha(\omega)$ is the largest open set with this property, we get $\mathscr{O}^k_\alpha\subseteq\tilde{Q}^k_\alpha(\omega)$.
\end{proof}

\begin{lemma}\label{lem:containedBall}
\begin{equation*}
  B(z^k_\alpha,\tfrac16 A_0^{-5}\delta^k)\subseteq\tilde{Q}^k_\alpha(\omega).
\end{equation*}
\end{lemma}

\begin{proof}
It suffices to prove that $B(z^k_\alpha,\tfrac16 A_0^{-5}\delta^k)$ is disjoint from each $\bar{Q}^k_\beta(\omega)$ with $\beta\neq\alpha$. For contradiction, let $x\in B(z^k_\alpha,\tfrac16 A_0^{-5}\delta^k)\cap\bar{Q}^k_\beta(\omega)$. Thus $x=\lim_{m\to\infty}z^m_{\theta(m)}$ for some $(m,\theta(m))\leq_\omega(k,\beta)$. Then
\begin{equation*}
  d(z^m_{\theta(m)},z^k_\alpha)
  \leq A_0 d(z^m_{\theta(m)},x)+A_0 d(x,z^k_\alpha)<\tfrac16A_0^{-4}\delta^k
\end{equation*}
for large enough $m$, since the second term is strictly smaller than this bound, and the first term tends to zero as $m\to\infty$. But then Lemma~\ref{lem:newLeqVsDistIter} says that $(m,\theta(m))\leq_\omega(k,\alpha)$, a contradiction with $(m,\theta(m))\leq_\omega(k,\beta)$.
\end{proof}

The following theorem summarizes the above properties of the dyadic cubes for a fixed parameter $\omega$, and supplements the key statement about their probabilistic behaviour under the random choice of $\omega\in\Omega$.

\begin{theorem}\label{thm:cubes}
For any fixed $\omega\in\Omega:=(\{0,1,\ldots,L\}\times\{1,\ldots,M\})^{\Z}$, the cubes satisfy the following relations: the covering properties
\begin{equation*}
  X=\bigcup_{\alpha}\bar{Q}^k_{\alpha}(\omega),\qquad
  \bar{Q}^k_{\alpha}(\omega)=\bigcup_{\beta:(k+1,\beta)\leq_{\omega}(k,\alpha)}\bar{Q}^{k+1}_{\beta}(\omega);
\end{equation*}
the mutual disjointness property
\begin{equation*}
  \bar{Q}^k_{\alpha}(\omega)\cap\tilde{Q}^k_{\beta}(\omega)=\varnothing\qquad(\alpha\neq\beta);
\end{equation*}
and the comparability with balls:
\begin{equation*}
  B(z^k_\alpha,\tfrac16 A_0^{-5}\delta^k)\subseteq\tilde{Q}^k_\alpha(\omega)
  \subseteq\bar{Q}^k_{\alpha}(\omega)\subseteq B(z^k_{\alpha},6A_0^4\delta^k).
\end{equation*}
Moreover, when $\Omega$ is equipped with the natural probability measure $\prob_{\omega}$, we have for some $\eta\in(0,1]$ the small boundary layer property:
\begin{equation}\label{eq:smallBdry}
  \prob_{\omega}\Big(x\in\bigcup_{\alpha}\partial_{\eps}Q^k_{\alpha}(\omega)\Big)\leq C\eps^{\eta}
  \qquad\Big(\partial_{\eps}Q^k_{\alpha}(\omega)
  :=\{y\in\bar{Q}^k_{\alpha}(\omega):d(y,{}^c\tilde{Q}^k_{\alpha}(\omega))<\eps\delta^k\}\Big);
\end{equation}
and in particular the negligible boundary property:
\begin{equation*}
  \prob_{\omega}\Big(x\in\bigcup_{k,\alpha}\partial Q^k_{\alpha}(\omega)\Big)=0
  \qquad\Big(\partial Q^k_{\alpha}(\omega):=\bar{Q}^k_{\alpha}(\omega)\setminus \tilde{Q}^k_{\alpha}(\omega)\Big).
\end{equation*}
\end{theorem}

\textbf{Henceforth, $\eta$ will always designate the fixed positive constant provided by this theorem.} It will reappear as the H\"older-regularity index of our splines and wavelets.

\begin{proof}
It only remains to check the probabilistic statements.

Let $x\in X$, $k\in\Z$, and $\varepsilon>0$ be fixed.
For every $\ell=k,k+1,\ldots$, there is some $\gamma$ so that $d(x,x^{\ell+1}_\gamma)<2A_0\delta^{\ell+1}$. Now suppose that for some $\ell$, the point $x^{\ell+1}_\gamma$ is chosen as a new dyadic point $z^\ell_\beta$, and recall that $B(z^\ell_\beta,\tfrac16 A_0^{-5}\delta^\ell)\subseteq\tilde{Q}^\ell_\beta(\omega)$. Thus,
\begin{equation*}
  d(x,{}^c\tilde{Q}^\ell_\beta(\omega))\geq
  \frac{1}{A_0}d(z^\ell_\beta,{}^c\tilde{Q}^\ell_\beta(\omega))-d(x,z^\ell_\beta)
  \geq\Big(\frac16 A_0^{-6}-2A_0\delta\Big)\delta^\ell
  \geq\frac17A_0^{-6}\delta^\ell
  \geq\varepsilon\delta^k
\end{equation*}
provided that $\ell\leq k+\log(7A_0^6\varepsilon)/\log\delta$. In particular, we have $x\in \tilde{Q}^\ell_\beta(\omega)\subseteq\tilde{Q}^k_\alpha(\omega)$ if $(\ell,\beta)\leq_\omega(k,\alpha)$, while $d(x,{}^c\tilde{Q}^k_\alpha(\omega))\geq d(x,{}^c\tilde{Q}^\ell_\beta(\omega))\geq\varepsilon\delta^k$. Thus $x$ \emph{cannot} be in $\bigcup_\alpha\partial_{\varepsilon}Q^k_\alpha(\omega)$ in this case.

In other words, in order that $x$ \emph{belongs} to the union $\bigcup_\alpha\partial_{\varepsilon}Q^k_\alpha(\omega)$, it is necessary that \emph{none} of the points $x^{\ell+1}_\gamma$, where $k\leq\ell\leq k+\log(7A_0^6\varepsilon)/\log\delta$, \emph{is chosen} as a new dyadic point. But, by Lemma~\ref{lem:prob}, every $x^{\ell+1}_\gamma$ has a probability at least $\tau:=\frac{1}{M(L+1)}$ of being chosen. So the probability that $x^{\ell+1}_\gamma$ is not chosen is at most $1-\tau$. Moreover, these events for different levels $\ell$ are independent of each other. Hence the probability that none of the $x^{\ell+1}_\gamma$ is chosen, for $k\leq\ell\leq k+\log(7A_0^6\varepsilon)/\log\delta$, is at most
\begin{equation*}
  (1-\tau)^{\log(7A_0^6\varepsilon)/\log\delta}
  =(7A_0^6\varepsilon)^{\log(1-\tau)/\log\delta}
  =C\varepsilon^{\eta},
\end{equation*}
where $\eta:=\log(1-\tau)/\log\delta>0$ (since both $\delta,1-\tau\in(0,1)$) and $C=(7A_0^6)^{\eta}$. This is exactly as claimed.

The negligible boundary property follows from the small boundary property as $\varepsilon\to 0$.
\end{proof}

We conclude this section with the observation that in our construction, the original dyadic point $x^k_\alpha$ may also be viewed as a `centre' of the new dyadic cubes $\bar{Q}^k_\alpha(\omega)$ for all $\omega\in\Omega$:

\begin{lemma}\label{lem:cubeVsBalls}
$B(x^k_{\alpha},\tfrac{1}{8}A_0^{-3}\delta^k)\subseteq\bar{Q}^k_{\alpha}(\omega)\subseteq\bar{B}(x^k_{\alpha},8A_0^5\delta^k)$ for any $\omega\in\Omega$.
\end{lemma}

\begin{proof}
By Lemma~\ref{lem:cubeVsCentre}, we have for any $x\in\bar{Q}^k_{\alpha}$ that
\begin{equation*}
  d(x,x^k_{\alpha})
  \leq A_0 d(x,z^k_{\alpha})+A_0 d(z^k_{\alpha},x^k_{\alpha})
  <6A_0^5\delta^k+2A_0^2\delta^k\leq 8A_0^5\delta^k,
\end{equation*}
which shows that $\bar{Q}^k_{\alpha}\subseteq B(x^k_{\alpha},8A_0^5\delta^k)$.

For the other inclusion, let $d(x,x^k_{\alpha})<\tfrac18 A_0^{-3}\delta^k$.
Since $X=\bigcup_{\eta}\bar{Q}^{k+1}_{\eta}(\omega)$, we have $x\in\bar{Q}^{k+1}_{\eta}(\omega)\subseteq B(x^{k+1}_{\eta},8A_0^5\delta^{k+1})$ for some $\eta$. But then
\begin{equation*}
  d(x^{k+1}_{\eta},x^k_{\alpha})
  \leq A_0 d(x^{k+1}_{\eta},x)+A_0 d(x,x^k_{\alpha})
  <8A_0^6\delta^{k+1}+\tfrac18A_0^{-2}\delta^k<\tfrac12A_0^{-1}\delta^k.  
\end{equation*}
By another application of \eqref{eq:leqProperties}, this implies that $(k+1,\eta)\leq(k,\alpha)$.

We want to prove that in fact $(k+1,\eta)\leq_{\omega}(k,\alpha)$, since then $x\in\bar{Q}^{k+1}_{\eta}(\omega)\subseteq\bar{Q}^k_{\alpha}(\omega)$. By \eqref{eq:newLeq}, the only potential obstacle to this is that $d(x^{k+1}_{\eta},z^k_{\gamma})<\tfrac14A_0^{-2}\delta^k$ for some $\gamma\neq\alpha$. So suppose this is the case, and recall that $z^k_{\gamma}=x^{k+1}_{\theta}$ for some $(k+1,\theta)\leq(k,\gamma)$. Then
\begin{equation*}
\begin{split}
  d(x^{k+1}_{\theta},x^k_{\alpha})
  &\leq A_0 d(z^k_{\gamma},x^{k+1}_{\eta})+A_0^2 d(x^{k+1}_{\eta},x)+A_0^2 d(x,x^k_{\alpha}) \\
  &<\tfrac14A_0^{-1}\delta^k+8A_0^7\delta^{k+1}+\tfrac{1}{8}A_0^{-1}\delta^k<\tfrac12A_0^{-1}\delta^k,
\end{split}
\end{equation*}
and hence $(k+1,\theta)\leq(k,\alpha)$ by \eqref{eq:leqProperties}. But also $(k+1,\theta)\leq(k,\gamma)$, hence $\gamma=\alpha$, a contradiction.

We have shown that an arbitrary $x\in B(x^k_{\alpha},\tfrac{1}{8}A_0^{-3}\delta^k)$ satisfies $x\in\bar{Q}^{k+1}_{\eta}(\omega)$ for some $(k+1,\eta)\leq_{\omega}(k,\alpha)$; thus $x\in\bar{Q}^k_{\alpha}(\omega)$.
\end{proof}

\section{Construction of splines}

The construction of splines on $X$, and the proof of their basic properties, is now an easy consequence of the preparations from the previous section.
For every $(k,\alpha)$, we define the \textbf{spline function}
\begin{equation}\label{eq:defSpline}
  s^k_{\alpha}(x):=\prob_{\omega}\Big(x\in\bar{Q}^k_{\alpha}(\omega)\Big).
\end{equation}

\begin{theorem}
The splines \eqref{eq:defSpline} satisfy the following properties: bounded support
\begin{equation}\label{eq:splineSupport}
  1_{B(x^k_{\alpha},\tfrac{1}{8}A_0^{-3}\delta^k)}(x)\leq s^k_{\alpha}(x)\leq 1_{B(x^k_{\alpha},8A_0^5\delta^k)}(x);
\end{equation}
the interpolation and reproducing properties
\begin{equation}\label{eq:splineSum}
  s^k_\alpha(x^k_\beta)=\delta_{\alpha\beta},\qquad
  \sum_{\alpha}s^k_{\alpha}(x)=1,\qquad
  s^k_{\alpha}(x)
  =\sum_{\beta}p^k_{\alpha\beta}\cdot s^{k+1}_{\beta}(x)
\end{equation}
where $\{p^k_{\alpha\beta}\}_\beta$ is a finitely nonzero set of nonnegative coefficients with $\sum_\beta p^k_{\alpha\beta}=1$;
and H\"older-continuity
\begin{equation*}
  \abs{s^k_{\alpha}(x)-s^k_{\alpha}(y)}
    \leq C\Big(\frac{d(x,y)}{\delta^k}\Big)^{\eta}.
\end{equation*}
\end{theorem}

Regarding \eqref{eq:splineSupport}, remark that it is a bit unusual that a spline be a non-zero constant  on part of its support.

\begin{proof}
The relations \eqref{eq:splineSupport} are immediate from Lemma~\ref{lem:cubeVsBalls}, and this implies in particular that $s^k_{\alpha}(x^k_{\alpha})=1$.

Since the boundaries have vanishing probability and $X=\bigcup_{\alpha}\bar{Q}^k_{\alpha}(\omega)$, it follows that
\begin{equation}\label{sum}
  \sum_{\alpha}s^k_{\alpha}(x)=\prob_{\omega}\Big(x\in X\Big)=1.
\end{equation}
Since these functions are nonnegative and $s^k_{\alpha}(x^k_{\alpha})=1$, it must be that $s^k_{\alpha}(x^k_{\beta})=0$ for $\beta\neq\alpha$, and hence in fact one has the \textbf{interpolation property}
\begin{equation}\label{interpolation}
  s^k_{\alpha}(x^k_{\beta})=\delta_{\alpha\beta}.
\end{equation}

From the identity
\begin{equation*}
   \bar{Q}^k_{\alpha}(\omega)=\bigcup_{\beta:(k+1,\beta)\leq_{\omega}(k,\alpha)}\bar{Q}^{k+1}_{\beta}(\omega),
\end{equation*}
(and using again the vanishing probability of the boundaries) we also have
\begin{equation*}
\begin{split}
  s^k_{\alpha}(x) &=\prob_{\omega}\Big(x\in\bigcup_{\beta:(k+1,\beta)\leq_{\omega}(k,\alpha)}\bar{Q}^{k+1}_{\beta}(\omega)\Big) \\
   &=\sum_{\beta}\prob_{\omega}\Big(\Big\{(k+1,\beta)\leq_{\omega}(k,\alpha)\Big\}\cap\Big\{x\in\bar{Q}^{k+1}_{\beta}(\omega)\Big\}\Big) \\
   &=\sum_{\beta}\prob_{\omega}\Big((k+1,\beta)\leq_{\omega}(k,\alpha)\Big)\prob_{\omega}\Big(x\in\bar{Q}^{k+1}_{\beta}(\omega)\Big) \\
   &=\sum_{\beta}\prob_{\omega}\Big((k+1,\beta)\leq_{\omega}(k,\alpha)\Big)s^{k+1}_{\beta}(x)=:\sum_{\beta}p^k_{\alpha\beta}\cdot s^{k+1}_{\beta}(x),
\end{split}
\end{equation*}
where the key third step used the independence of the two events; namely, the event $(k+1,\beta)\leq_{\omega}(k,\alpha)$ depends only on $\omega_k$, while the cube $\bar{Q}^{k+1}_{\beta}(\omega)$ depends on $\omega_{\ell}$ for $\ell\geq k+1$. The support properties of the splines readily imply that only boundedly many of the coefficients $p^k_{\alpha\beta}$ are nonzero for a given $(k,\alpha)$, so that $\lspan\{s^k_{\alpha}\}_{\alpha}\subseteq\lspan\{s^{k+1}_{\beta}\}_{\beta}$.

The \textbf{H\"older-continuity of the splines} follows from the probabilistic smallness of the boundary regions, as expressed by \eqref{eq:smallBdry}. Indeed,
\begin{equation*}
\begin{split}
  \abs{s^k_{\alpha}(x)-s^k_{\alpha}(y)}
  &=\Babs{\prob_{\omega}\Big(x\in\bar{Q}^k_{\alpha}(\omega)\Big)-\prob_{\omega}\Big(y\in\bar{Q}^k_{\alpha}(\omega)\Big)} \\
  &=\Babs{\int_{\Omega}\big(1_{\omega:x\in \bar{Q}^k_{\alpha}(\omega)}-1_{\omega:y\in \bar{Q}^k_{\alpha}(\omega)}\big)\ud\prob_{\omega}} \\
  &\leq \int_{\Omega}\big(1_{\omega:x\in \bar{Q}^k_{\alpha}(\omega),y\notin\bar{Q}^k_{\alpha}(\omega) }
     +1_{\omega:y\in \bar{Q}^k_{\alpha}(\omega),x\notin\bar{Q}^k_{\alpha}(\omega) }\big)\ud\prob_{\omega} \\
  &= \prob_{\omega}\Big(x\in \bar{Q}^k_{\alpha}(\omega),y\notin\bar{Q}^k_{\alpha}(\omega)\Big)
     +\prob_{\omega}\Big(y\in \bar{Q}^k_{\alpha}(\omega),x\notin\bar{Q}^k_{\alpha}(\omega) \Big) \\
  &\leq \prob_{\omega}\Big(x\in \partial_{d(x,y)\delta^{-k}}\bar{Q}^k_{\alpha}(\omega)\Big)
     +\prob_{\omega}\Big(y\in \partial_{d(x,y)\delta^{-k}}\bar{Q}^k_{\alpha}(\omega)\Big) \\
    &\leq C\Big(\frac{d(x,y)}{\delta^k}\Big)^{\eta}.\qedhere
\end{split}
\end{equation*}
\end{proof}

{
\section{Auxiliary results on quasi-metric spaces}
}

We include this intermediate section to collect some auxiliary material so as to streamline the subsequent analysis. This section is mostly concerned with difficulties of a quasi-metric in contrast to a metric. The results are mostly part of the folklore, but somewhat difficult to find in the literature in full generality, since additional assumptions on the space are usually imposed. { (After completion of this work, we became aware of  the recent preprint \cite{MMMM} where the density result below is also proved in full generality by a completely different argument.)}

For the beginning of this section, as above, we only assume that $X$ is a geometrically doubling quasi-metric space. \textbf{Throughout, let $\eta>0$ denote the fixed positive constant from the small boundary property of the dyadic cubes and the H\"older regularity of the splines}.

\begin{lemma}\label{lem:existsLip}
Let $F\subseteq G\subseteq X$ be sets with
\begin{equation*}
  d(F,G^c):=\inf_{x\in F,y\notin G}d(x,y)=\Delta>0.
\end{equation*}
Then there exists a function $\varphi:X\to\R$ with $1_{F}\leq\varphi\leq 1_G$ and
\begin{equation*}
  \abs{\varphi(x)-\varphi(y)}\leq C\Big(\frac{d(x,y)}{\Delta}\Big)^{\eta}.
\end{equation*}
\end{lemma}

\begin{proof}
Let $k$ be the smallest integer so that $16A_0^6\delta^k\leq\Delta$. We set
\begin{equation*}
  \varphi:=\sum_{\alpha:B(x^k_\alpha,8A_0^5\delta^k)\cap F\neq\varnothing}s^k_\alpha.
\end{equation*}
Then $\varphi(x)=1$ for $x\in F$, since $\sum_\alpha s^k_\alpha(x)=1$, and the sum defining $\varphi(x)$ contains all $s^k_\alpha(x)$ whose support intersects $F$. To prove that $\varphi(z)=0$ for $z\in G^c$, we need to check that no $\supp s^k_\alpha\subseteq B(x^k_\alpha,8A_0^5\delta^k)$ appearing in $\varphi$ can intersect $G^c$. To this end, let $y\in B(x^k_\alpha,8A_0^5\delta^k)\cap F$ (which exists by definition) and $z\in B(x^k_\alpha,8A_0^5\delta^k)\cap G^c$ (whose assumed existence should lead to a contradiction). Then
\begin{equation*}
  d(F,G^c)\leq d(y,z)< A_0 d(y,x^k_\alpha)+A_0 d(x^k_\alpha,z)\leq 2A_0\cdot 8 A_0^5\delta^k\leq\Delta,
\end{equation*}
a contradiction indeed. So we have $1_F\leq\varphi\leq 1_G$, and it remains to check the H\"older-continuity. By geometric doubling, there are only boundedly many $s^k_\alpha$, whose support contains either $x$ or $y$. Each of these $s^k_\alpha$ satisfy the estimate
\begin{equation*}
  \abs{s^k_\alpha(x)-s^k_\alpha(y)}\leq C\Big(\frac{d(x,y)}{\delta^k}\Big)^{\eta}\leq C\Big(\frac{d(x,y)}{\Delta}\Big)^{\eta},
\end{equation*}
and hence so does their sum over boundedly many indices $\alpha$. This completes the proof.
\end{proof}

\begin{corollary}\label{cor:existBump}
For every $B(x,r)$, there exists a function $\varphi:X\to\R$ with $1_{B(x,r)}\leq\varphi\leq 1_{B(x,2A_0 r)}$ and
\begin{equation*}
    \abs{\varphi(x)-\varphi(y)}\leq C\Big(\frac{d(x,y)}{r}\Big)^{\eta}.
\end{equation*}
\end{corollary}

\begin{proof}
If $y\in B(x,r)$ and $z\in B(x,2A_0 r)^c$, then
\begin{equation*}
  d(y,z)\geq\frac{1}{A_0}d(z,x)-d(x,y)>2r-r=r;
\end{equation*}
thus $d(B(x,r),B(x,2A_0r)^c)\geq r$, and the previous lemma applies.
\end{proof}

We formulate a quasi-metric version of a well-known covering lemma for metric spaces, cf. \cite[Theorem 1.2]{Heinonen}; the quasi-metric extension is obtained \emph{mutatis mutandis}.

\begin{lemma}
Let $\mathscr{B}$ be a family of balls $B(x,r)$ in $X$, with bounded radii and contained in a bounded set. Then there exists a pairwise disjoint subcollection $\{B_i\}_{i=1}^{\infty}$, whose concentric expansions $5A_0^2 B_i$ cover all original $B(x,r)\in\mathscr{B}$.
\end{lemma}

\begin{lemma}\label{lemma:union}
$X$ is a countable union of bounded open sets.
\end{lemma}

\begin{proof}
Recall that the interior of a set $E$ is the open set $(\overline{E^c})^c\subseteq E$. It is easy to check that
\begin{equation*}
  B(x,r/A_0)\subseteq\operatorname{interior}B(x,r);
\end{equation*}
hence, with any base point $x_0\in X$, the bounded open sets $\operatorname{interior}B(x_0,n)$, $n\in\N$, provide the required covering.
\end{proof}

So far everything has been based on the space geometry, i.e., the properties of the quasi-distance only and geometric doubling. We next add a {measure $\mu$} into our considerations.  

\begin{proposition}\label{prop:density} Let $\mu$ be a non-trivial Borel measure on $X$, finite on bounded Borel sets.
Let $1\leq p<\infty$. Then H\"older-$\eta$-continuous functions of bounded support are dense in $L^p(\mu)$, where $\eta$ is the H\"older exponent of the splines.
\end{proposition}

\begin{proof}
By the density of simple functions due to general measure theory, it suffices to show that for every bounded Borel set $E$ and every $\epsilon>0$, there exists a boundedly supported H\"older-$\eta$-continuous function $\varphi$ with $\Norm{1_E-\varphi}{p}<2\epsilon$. By a general result concerning Borel measures \cite[Theorem 2.2.2 (ii)]{Federer}\footnote{This result is stated in a metric space with an outer measure but the argument is valid for any  topological space for a Borel measure. Moreover, it is assumed that $E$ can be covered by countably many open sets with bounded measure which is granted from Lemma \ref{lemma:union} and our assumptions on $\mu$.}, there is an open set $G\supseteq E$ 
such that $\mu(G\setminus E)<\epsilon^p$. 
For every $x\in G$, we choose a radius $r_x$ small enough so that $B(x,4A_0^2 r_x)\subseteq G$. The balls $B(x,5^{-1}A_0^{-2}r_x)$ cover $G$, and hence there are pairwise disjoint balls $B(x_i,5^{-1}A_0^{-2}r_i)$, $r_{i}:=r_{x_{i}}$,  so that $B(x_{i},r_{i})$ cover $G$.  As we have $B(x,r) \subseteq \overline{B(x,r)} \subseteq B(x,2A_{0}r)$ for any ball, with $\overline B$ being the topological closure of the ball $B$,
\begin{equation*}
  G\subseteq\bigcup_{i=1}^{\infty}\overline{B(x_i,r_i)}\subseteq\bigcup_{i=1}^{\infty} B(x_i,2A_0 r_i)\subseteq G.
\end{equation*}
 Thus
 \begin{equation*}
 \mu(G)=\lim_{n\to\infty}\mu\Big(\bigcup_{i=1}^n \overline{B(x_i,r_i)}\Big)=:\lim_{n\to\infty}\mu(G_n).
\end{equation*}
Let us fix $n$ so large that $\mu(G\setminus G_n)\leq\epsilon^p$. We check that $d(G_n,G^c)>0$. Indeed, if $x\in G_n$, then $x\in\overline{B(x_i, r_i)}\subseteq B(x_i,2A_0 r_i)$ for some $i=1,\ldots,n$, whereas $B(x_i,4A_0^2r_i)\subseteq G$. Hence
\begin{equation*}
  d(G_n,G^c)\geq\min_{i=1,\ldots,n} d(B(x_i,2A_0r_i),B(x_i,4A_0^2 r_i)^c)\geq\min_{i=1,\ldots,n} { r_i}>0.
\end{equation*}
By Lemma~\ref{lem:existsLip}, there exists a H\"older-$\eta$-continuous $\varphi$ with $1_{G_n}\leq\varphi\leq 1_G$. Hence $\Norm{1_G-\varphi}{p}\leq\Norm{1_G-1_{G_n}}{p}=\mu(G\setminus G_n)^{1/p}<\epsilon$, and finally $\|1_{E}-\varphi\|_{p}<2\epsilon$.
\end{proof}

\section{$L^2$ theory and multiresolution analysis}\label{sec:l2}

We now return to the development of the spline theory in the presence of a nontrivial Borel measure $\mu$ on $(X,d)$. From now on, we assume that $(X,d,\mu)$ is a space of homogeneous type in the most general sense, namely that in addition to having a quasi-metric $d$ we only assume the doubling condition    which reads:   for all $x\in X$ and $r>0$ 
\begin{equation*}
0<  \mu(B(x,2r))\leq C_{\mu}\cdot\mu(B(x,r))<\infty.
\end{equation*}
This inequality makes sense if balls are Borel sets. But this may not be the case.  However, their topological closures are (because they are closed sets by definition) and satisfy $B(x,r) \subseteq \overline{B(x,r)} \subseteq B(x,2A_{0}r)$. In that case, up to changing constants throughout, we may change balls to their closures in the above definition of the doubling condition, thus avoiding to resort to the outer measure associated to $\mu$. To simplify matters though, we assume that balls are Borel sets and leave to the reader the modifications if not.

The doubling condition on $(X,d,\mu)$ implies that $(X,d)$ is geometrically doubling (see \cite{CW}), so that  the spline functions, built independently of measure considerations, are  at our disposal. 
We now show that they provide a multiresolution analysis of $L^2(\mu)$, {essentially in the sense of Sweldens \cite[Definition 3.1]{Sweldens}, who considered the case of general measures on $\R^n$. This consist of all properties of a classical multiresolution analysis of Meyer \cite[Definition 2.1]{M2}, to the extent that this definition is meaningful in a quasi-metric space context: the classical postulates dealing with translations and dilations, specific to the Euclidean space and the Lebesgue measure, are now meaningless.}

\begin{theorem}\label{thm:multires}
Let $V_k$ be the closed linear span of $\{s^k_\alpha\}_\alpha$ {in $L^2(\mu)$}. Then $V_k\subseteq V_{k+1}$, and
\begin{equation*}
  \overline{\bigcup_{k\in\Z}V_k}=L^2(\mu),\qquad
  \bigcap_{k\in\Z}V_k=\begin{cases} \{0\}, & \text{if $X$ is unbounded}, \\ V_{k_0}=\{\mathrm{constants}\}, & \text{if $X$ is bounded}, \end{cases}
\end{equation*}
where $k_0$ is some integer.
Moreover, the functions $s^k_\alpha/\sqrt{\mu^k_{\alpha}}$ form a Riesz basis of $V_k$: for all sequences of numbers $\lambda_{\alpha}$, we have the two-sided estimate
\begin{equation*}
  \BNorm{\sum_{\alpha}\lambda_{\alpha}s^k_{\alpha}}{L^2(\mu)}\eqsim\Big(\sum_{\alpha}\abs{\lambda_{\alpha}}^2 \mu^k_{\alpha}\Big)^{1/2},
\end{equation*}
with $\mu^k_{\alpha}:= \mu(B(x^k_{\alpha},\delta^k)).$
\end{theorem}

\begin{proof}
The nesting property $V_{k}\subseteq V_{k+1}$ is immediate from the reproducing properties of the splines.

\subsubsection*{The Riesz basis property}
We use properties \eqref{eq:splineSupport} and \eqref{eq:splineSum} of the splines.
At any point $x$, the values $s^k_{\alpha}(x)$ are nonnegative and sum up to $1$ (when $\alpha$ is the summation variable and $k$ is fixed). Hence
\begin{equation*}
  \Babs{\sum_{\alpha}\lambda_{\alpha}s^k_{\alpha}(x)}^2
  \leq\sum_{\alpha}\abs{\lambda_{\alpha}}^2 s^k_{\alpha}(x)\leq \sum_{\alpha}\abs{\lambda_{\alpha}}^2  1_{B(x^k_{\alpha},8A_0^5\delta^k)}(x).
\end{equation*}
Thus
\begin{equation*}
  \int\Babs{\sum_{\alpha}\lambda_{\alpha}s^k_{\alpha}(x)}^2\ud\mu(x)
  \leq\sum_{\alpha}\abs{\lambda_{\alpha}}^2  \mu(B(x^k_{\alpha},8A_0^5\delta^k))
  \leq C\sum_{\alpha}\abs{\lambda_{\alpha}}^2  \mu(B(x^k_{\alpha},\delta^k))
\end{equation*}
by the doubling property.

On the other hand, we also know that $s^k_{\beta}(x)=1$ for $x\in B(x^k_{\beta},\tfrac{1}{8A_{0}^3}\delta^k)$ and all other $s^k_{\alpha}$, $\alpha\neq\beta$, vanish on this set. Hence
\begin{equation*}
  \Babs{\sum_{\alpha}\lambda_{\alpha}s^k_{\alpha}(x)}^2
  =\abs{\lambda_{\beta}}^2,\qquad\forall x\in B(x^k_{\beta},\tfrac{1}{8A_{0}^3}\delta^k),
\end{equation*}
thus
\begin{equation*}
  \Babs{\sum_{\alpha}\lambda_{\alpha}s^k_{\alpha}(x)}^2
  \geq\sum_{\alpha}\abs{\lambda_{\alpha}}^2 1_{B(x^k_{\alpha},\tfrac{1}{8A_{0}^3}\delta^k)},
\end{equation*}
and therefore
\begin{equation*}
  \int \Babs{\sum_{\alpha}\lambda_{\alpha}s^k_{\alpha}(x)}^2\ud\mu(x)
  \geq\sum_{\alpha}\abs{\lambda_{\alpha}}^2 \mu(B(x^k_{\alpha},\tfrac{1}{8A_{0}^3}\delta^k))
  \geq\frac{1}{C}\sum_{\alpha}\abs{\lambda_{\alpha}}^2 \mu(B(x^k_{\alpha},\delta^k)),
\end{equation*}
again by the doubling property.

\subsubsection*{The union of the  $V_{k}$}
If  $f$ is a H\"older-$\eta$-continuous function with bounded support, the sum 
$f_{k}(x)= \sum_{\alpha}f(x^k_{\alpha}) s^k_{\alpha}(x)$
defines an element of $V_{k}$ and using \eqref{eq:splineSum} we have
$$
f(x)-f_{k}(x)= \sum_{\alpha}(f(x)-f(x^k_{\alpha}) )s^k_{\alpha}(x)
$$
so that 
$$
\|f-f_{k}\|_{\infty} \le \sup_{\alpha}\sup_{x\in B(x^k_{\alpha}, 8A_{0}^6\delta ^k)}|f(x)-f(x^k_{\alpha})| \leq C \delta ^{k\eta}.
$$
 Convergence  in  $L^2(\mu)$ follows from this and bounded support for $f-f_{k}$.

\subsubsection*{The intersection of the $V_k$}
In case $X$ is bounded, it is clear from definition that the sets $\{x^k_{\alpha}\}_{\alpha}$ reduce to one point when $k$ gets small. In that case, $s^k_{\alpha}(x)=1$ for all $x\in X$ by \eqref{eq:splineSum}. So we choose $k_{0}$ as the largest 
integer  such that this property holds.
 
If $X$ is unbounded,  then $\mu(X)=\infty$ as $\mu$ is doubling. (See remark after Lemma~\ref{lem:emptyAnnulus} for a proof of this known fact.) If $f \in  \bigcap_{k\in \Z} V_{k}$, then 
 $\|f\|_{2}^2 \sim  \sum_{\alpha} |f(x^k_{\alpha})|^2\mu^k_{\alpha}$  for all $k$.
 Consider a fixed $x\in X$. Then $f(x)=\sum_\alpha f(x^k_\alpha)s^k_\alpha(x)$, where the sum is over the boundedly (with respect to $k$) many $\alpha$ such that $x^k_\alpha\in B(x,8A_0^5\delta^k)$. Thus $\abs{f(x)}\leq\sum_\alpha\abs{f(x^k_\alpha)}\lesssim\sum_\alpha\Norm{f}{2}/\sqrt{\mu^k_\alpha}\lesssim\Norm{f}{2}/V(x,\delta^k)$. This tends to zero as $k\to-\infty$, and hence $f(x)=0$. Since the point was arbitrary, we have $f\equiv 0$.
 \end{proof}

\section{Biorthogonal and orthogonal spline systems}\label{sect:biortho}

In this section, we use a classical algorithm (cf. \cite[Sec.~2.3]{M2}) to construct two further bases of the space $V_k$ spanned by the splines $\{s^k_\alpha\}_\alpha$. We use the abbreviations
\begin{equation*}
  \mu^k_\alpha:=V(x^k_\alpha,\delta^k):=\mu(B(x^k_\alpha,\delta^k))
\end{equation*}
for these frequently appearing volumes.

\begin{theorem}\label{thm:biortho}
There exist a system of biorthogonal splines $\{\tilde{s}^k_\alpha\}_\alpha$ in $V_k$ with
\begin{equation}\label{eq:biortho}
\langle s^k_{\alpha}, \tilde s^k_{\beta}\rangle_{L^2(\mu)}= \delta _{\alpha,\beta},
\end{equation}
as well as an orthonormal basis $\{\phi^k_\alpha\}_\alpha$ of $V_k$, which satisfy the following estimates, 
where 
$d(x,y)\leq\delta^k$:
\begin{equation*}
\begin{split}
  \mu^k_{\alpha}\cdot\abs{\tilde s^k_{\alpha}(x)}
  +\sqrt{\mu^k_{\alpha}}\cdot\abs{\phi^k_{\alpha}(x)}
  &\leq C\exp\big(-\gamma {{\delta^{-k}}{d(x^k_{\alpha},x)}}^s\, \big) \\
  \mu^k_{\alpha}\cdot\abs{\tilde s^k_{\alpha}(x)-\tilde s^k_\alpha(y)}
  +\sqrt{\mu^k_{\alpha}}\cdot\abs{\phi^k_{\alpha}(x)-\phi^k_\alpha(y)}
  &\leq C\Big(\frac{d(x,y)}{\delta^k}\Big)^{\eta} \exp\big( -\gamma (\delta ^{-k}{d(x^k_{\alpha},x)})^s\, \big),
\end{split}
\end{equation*}
with  $s=(1+\log_{2}A_{0})^{-1}$ or $s=1$ if $d$ is a Lipschitz-continuous quasi-distance.
\end{theorem}

\begin{remark}
From the decay estimates and doubling, it readily follows that $\tilde{s}^k_\beta\in L^1(\mu)$. Summing the biorthogonality relation~\eqref{eq:biortho} over all $\alpha$ and recalling that $\sum_\alpha s^k_\alpha(x)\equiv 1$, we deduce that
\begin{equation*}
  \int \tilde{s}^k_\beta\ud\mu = 1.
\end{equation*}
\end{remark}

We begin the proof of Theorem~\ref{thm:biortho}. Fix $k\in\Z$. By abuse, we identify $\mathscr{X}^k=\{x^k_{\alpha}\}_{\alpha}$ and the set of indices $\alpha$ corresponding to points $x^k_{\alpha}$. The previous section says that there is a linear,  bounded, injective map
$U_{k}\colon\ell^2(\mathscr{X}^{k}) \to L^2(\mu)$  with closed range, defined by \begin{equation*}U_{k}\lambda= \sum_{\alpha} \frac{\lambda_{\alpha}}{\sqrt {\mu^k_{\alpha}}} s^k_{\alpha}
\end{equation*} for 
$\lambda=\{\lambda_{\alpha}\}_{\alpha\in \mathscr{X}^{k}}$.  Here, $\ell^2(\mathscr{X}^{k})$  is equipped with counting measure. Let $V_{k}$ denote its  range. 
By the properties of the splines, the sum defining each element of $V_{k}$ converges locally uniformly and  defines a H\"older-$\eta$-continuous function. 
The inverse of $U_{k}$ can be computed using \eqref{interpolation} by 
$(U_{k}^{-1}f)_{\alpha}= f(x^k_{\alpha}) \sqrt{\mu^k_{\alpha}}$. 

Let $\delta ^k_{\alpha}$ be the canonical orthonormal basis element  in $\ell^2(\mathscr{X}^{k})$. By construction  $U_{k}\delta ^k_{\alpha}= \frac{s^k_{\alpha}}{\sqrt {\mu^k_{\alpha}}}$. Since $U_{k}$ is an isomorphism, this means that the splines $s^k_{\alpha}$ form an unconditional basis of $V_{k}$. To find 
a biorthogonal system $\tilde s^k_{\alpha}$ in $V_{k}$, that is a system such that 
$$
\langle s^k_{\alpha}, \tilde s^k_{\beta}\rangle_{L^2(\mu)}= \delta _{\alpha\beta},
$$ we observe that if 
 $f=U_{k}\lambda, f'=U_{k}\lambda'$ then
\[
\langle f, f'\rangle_{L^2(\mu)}= \langle M_{k}\lambda, \lambda'\rangle_{\ell^2(\mathscr{X}^{k})}
\]
with $M_{k}$ being the infinite matrix with entries
$$
M_{k}(\alpha,\beta)=\frac{\langle s^k_{\alpha}, s^k_{\beta}\rangle_{L^2(\mu)}}{\sqrt{\mu^k_{\alpha}\mu^k_{\beta}}}.
$$
Another reformulation of the isomorphism property of $U_{k}$ is that $M_{k}$ is bounded and invertible on $\ell^2(\mathscr{X}^{k})$. It is also positive self-adjoint. So the biorthogonal system is uniquely defined by
$$
\sqrt{\mu^k_{\alpha}}\tilde s^k_{\alpha} = U_{k}M_{k}^{-1}\delta ^k_{\alpha}.
$$
If one wants an orthonormal basis of $V_{k} \subset L^2(\mu)$, one defines instead
$$
\phi^k_{\alpha}=U_{k}M_{k}^{-1/2} \delta ^k_{\alpha}.
$$
In other words,
\begin{equation*}
\begin{split}
\sqrt{\mu^k_{\alpha}}\tilde s^k_{\alpha}(x) &= \sum_{\beta\in \mathscr{X}^{k}} M_{k}^{-1} (\alpha,\beta)  \frac{s^k_{\beta}(x)}{\sqrt{\mu^k_{\beta}}}, \\
\phi^k_{\alpha}(x) &= \sum_{\beta \in \mathscr{X}^{k}} M_{k}^{-1/2} (\alpha,\beta)  \frac{s^k_{\beta}(x)}{\sqrt{\mu^k_{\beta}}}.
\end{split}
\end{equation*} 

We now establish decay estimate of the coefficients $M_k^{-1}(\alpha,\beta)$ and $M_k^{-1/2}(\alpha,\beta)$.

\begin{proposition}\label{prop:matrices}
There exist constants $\gamma>0$ and $C<\infty$, independent of $k$ such that the entries of $M_{k}^{-1}$ and $M_{k}^{-1/2}$ have upper bounds
$$
C\exp\big(- \gamma (\delta ^{-k}{d(\alpha,\beta)})^{s}\, \big),
$$
with $s=(1+\log_{2}A_{0})^{-1}$ {or $s=1$ if $d$ is Lipschitz-continuous.}
\end{proposition}

If we introduce the induced normalized quasi-distance on $\mathscr{X}^{k}\times \mathscr{X}^{k} $,   $d_{k}(\alpha,\beta):= \frac{d(x^k_{\alpha},x^k_{\beta})}{\delta ^k}$, we have to prove uniform  estimates on the entries of $M_{k}^{-1}$ or $M_{k}^{-1/2}$ in terms of $d_{k}$.  Note that for this quasi-distance $M_{k}$ is a band-matrix, more precisely $M_{k}(\alpha,\beta)= 0$ if $d_{k}(\alpha,\beta) \ge 16A_{0}^6 $ by \eqref{eq:splineSupport}, hence it has exponential decay.  

\begin{lemma}\label{lemma1} Let $\Xi$ be a $1$-separated set in   a quasi-metric space $(X,d)$ with quasi-triangle constant $A_{0}$ having the geometric doubling property with constant $N$. Then   for all $\varepsilon>0$, there exists $c(\varepsilon,A_{0},N)<\infty$ such that 
$$
\sup_{\alpha \in X}{\exp\big( \varepsilon d(\alpha, \Xi)/A_{0}\big)} \sum_{\beta \in \Xi} \exp\big(- \varepsilon d(\alpha,\beta)\big) \le c(\varepsilon,A_{0},N). 
$$
 Moreover,
for any $0<c'<c/A_{0}$
$$
 \sum_{\beta \in \Xi}  \exp\big(- c d(\alpha,\beta)\big) \exp\big(- c d(\beta,\gamma)\big) \le C \exp\big(- c' d(\alpha,\gamma)\big)
$$
where $C$ does not depend on $\alpha,\gamma \in X$.
\end{lemma}

\begin{proof} Fix a point $\alpha \in X$.  Pick any point $\beta'\in \Xi$. We have
$$
 \sum_{\beta \in \Xi} \exp\big(- \varepsilon d(\alpha,\beta)\big) \le  \exp\big(- \varepsilon d(\alpha,\beta')/A_{0}\big)   \sum_{\beta \in \Xi} \exp\big(- \varepsilon d(\beta',\beta)/A_{0}\big).
 $$
The number  of points in $\Xi$ at quasi-distance 
at most $2^j$ from $\beta'$ is bounded by $C^{j+1}$ from some $C$ depending only on $A_{0}$ and $N$ but not $\beta'$. Thus
$$\sum_{\beta \in \Xi} \exp\big(- \varepsilon d(\beta',\beta)/A_{0}\big) \le C+  \sum_{j\ge 1}  \sum_{\beta: d(\beta',\beta)\sim 2^j}  \exp\big(- \varepsilon d(\beta',\beta)/A_{0}\, \big)  \le \sum_{j\ge 0} C^{j+1} \cdot
\exp\big(- \varepsilon 2^j/A_{0}\, \big).$$
Taking the infimum over all $\beta'\in \Xi$ proves the desired inequality.

The result for the matrix product coefficients follows analogously and we skip details. 
\end{proof}

\begin{lemma}\label{lemma2} Let $\Xi$ be a $1$-separated set in   a quasi-metric space $(X,d)$ with quasi-triangle constant $A_{0}$ having the geometric doubling property with constant $N$. Consider a matrix  $M=(M(\alpha,\beta))$ indexed by $\Xi\times \Xi$ such that there exists $c>0$ for which
$$
C=\sup_{(\alpha,\beta)} \exp\big(c d(\alpha,\beta)\big) |M(\alpha,\beta)| <\infty.
$$
Then $M$ is bounded on $\ell^2(\Xi)$. 
If $M$ is invertible, then there exists $c'>0$ such that  
$$
\sup_{(\alpha,\beta)} \exp\big(c' (d(\alpha,\beta))^s\, \big) |M^{-1}(\alpha,\beta)| <\infty
$$
with  $s=(1+\log_{2} A_{0})^{-1}$ or $s=1$ if $d$ is a Lipschitz-continuous quasi-distance. If, in addition, $M$ is positive self-adjoint, then the same conclusion holds for $M^{-1/2}$.
\end{lemma}

\begin{remark}
(i) Note that the exponent $s=1$ for a usual distance is recovered as a special case in two ways, either by setting $A_{0}=1$ or using the Lipschitz-continuity.

{
(ii) The exponent $s$ is in general optimal in this result. Namely, consider the band matrix $M$ indexed by $\Z$ with $M(i,i)=1$ and $M(i,i+1)=-\lambda\in(-1,0)$ for all $i\in\Z$, and all other entries equal to zero. Then $M^{-1}(i,j)=\lambda^{i-j}$ if $j\geq i$ and zero otherwise. Now equip $\Z$ with the quasi-distance $d(i,j)=\abs{i-j}^r$, where $r\geq 1$. Its quasi-triangle constant is $A_0=2^{r-1}$, and hence the exponent given by the lemma is $s=(1+\log_2 A_0)^{-1}=(1+r-1)^{-1}=1/r$. This gives precisely the correct decay of the matrix with respect to $d$, as
\begin{equation*}
  \lambda^{i-j}=\exp(-c\cdot\abs{i-j})=\exp(-c\cdot d(i,j)^{1/r}),\qquad c=\log\lambda^{-1}.
\end{equation*}
}
\end{remark}

\begin{proof}  If $d$ is a genuine distance or a Lipschitz continuous quasi-distance, then this follows from Theorem 5 in \cite{Lem2} and the remark that follows it, which extends an earlier result in \cite{Dem} for band-limited matrices (but we shall need the full result here), with $s=1$.  
However, as said not all quasi-distances are Lipschitz- or even H\"older-continuous and we provide an argument {in full generality, which also recovers the mentioned special cases.}

 We begin with the following observation. For $n\ge 1$, let $\kappa_{n}$ be the best  constant in the inequality
$$
d(\alpha_{0},\alpha_{n}) \le \kappa_{n }(d(\alpha_{0},\alpha_{1})+ d(\alpha_{1},\alpha_{2})+ \ldots +d(\alpha_{n-1},\alpha_{n}))$$
for every chain $(\alpha_{0},\alpha_{1}, \dots, \alpha_{n})$ of $n+1$ elements (not necessarily distinct) of $\Xi$. It is clear that $(\kappa_{n})$ is non-decreasing and $\kappa_{1}=1$, $\kappa_{2}\le A_{0}$. 
Moreover, using $d(\alpha_{0}, \alpha_{m+n}) \le A_{0}(d(\alpha_{0}, \alpha_{m})+ d(\alpha_{m}, \alpha_{m+n}))$, it follows that $\kappa_{m+n}\le A_{0}{\cdot\max\{\kappa_{m},\kappa_{n}\}}$. {Thus $\kappa_{2n}\le A_{0}\kappa_{n}$ and therefore, $\kappa_{2^j}\le A_{0}^j$. 
We conclude that $\kappa_{n}\le A_{0}^{1+\log_{2} n}= A_{0} n ^{\log_{2}A_{0}}$. Note also that if $d$ is $L$-Lipshitz, then $d(\alpha_0,\alpha_n)\leq d(\alpha_0,\alpha_{n-1})+Ld(\alpha_{n-1},\alpha_n)$, which by iteration gives $\kappa_n\leq L$ for all $n$.}

Now assume that $M$ is {positive} self-adjoint and invertible. In this case, one can write $M=h(I-A)$ with $h=(\|M\| + \|M^{-1}\|)/2$ a positive real number and $A$ a matrix with norm $r=(\|M\|- \|M^{-1}\|)/(\|M\| + \|M^{-1}\|)<1$. Moreover, the coefficients of $A$ have the same decay as those of $M$. Without loss of generality, we normalize $h=1$.  Develop 
$
(I-A)^{-1}$ in the Neumann series $\sum A^n$ and estimate the coefficients $A^n(\alpha,\beta)$, $n\ge 1$, $\alpha\ne \beta$, in two ways. First 
$|A^n(\alpha,\beta)|\le r^n$. Second, we have {
\begin{equation*}
\begin{split}
  |A^n(\alpha,\beta)|
  &\leq \sum_{(\alpha_{1}, \ldots, \alpha_{n-1})\in \Xi^{n-1}} C^n \exp \big(-c (d(\alpha,\alpha_{1})+ d(\alpha_{1},\alpha_{2})+ \ldots +d(\alpha_{n-1},\beta))\big) \\ 
  &\leq C^n\exp\big(-\frac{c}{2\kappa_n}d(\alpha,\beta)\big)\sum_{(\alpha_{1}, \ldots, \alpha_{n-1})\in \Xi^{n-1}} 
   \exp \big(-\frac{c}{2} (d(\alpha_{1},\alpha_{2})+ \ldots +d(\alpha_{n-1},\beta))\big) \\
  &\leq \tilde{C}^n\exp\big(-\frac{c}{2\kappa_n}d(\alpha,\beta)\big),
\end{split}
\end{equation*}
where we applied Lemma~\ref{lemma1} $n-1$ times with $\varepsilon=c/2$, and $\tilde{C}=C\cdot c(\varepsilon,A_0,N)$, in the notation of that lemma.
 }
As $\kappa_{n}$ is non decreasing, we have for any  integer $n_{0}$ using the second estimate for $0\le n\le n_{0}$ and the first for $n>n_{0}$,
\begin{align*}
 |M^{-1}(\alpha,\beta)| &\le (n_{0}+1) \tilde C^{n_{0}} \exp\big(-\frac{c}{2\kappa_{n_{0}}} d(\alpha,\beta)\big) + r^{n_{0}+1} (1-r)^{-1} 
\end{align*}
{and $(n_{0}+1) \tilde C^{n_{0}}\leq D^{n_0}$} for some large constant $D>0$.
Choosing $n_{0}$ as the first integer such that the first term dominates, we see that $d(\alpha,\beta)\eqsim n_0\cdot\kappa_{n_0}\lesssim n_0^{1/s}$, where $s=1/(1+\log_2 A_0)$ (or $s=1$ if $d$ is Lipschitz-continuous, recalling that $\kappa_{n_0}\leq L$ in this case). Hence, for some constant $c'>0$, 
\begin{equation*}
  \abs{M^{-1}(\alpha,\beta)}\lesssim r^{n_0}
  \lesssim\exp\big(-c' d(\alpha,\beta)^s\, \big).
\end{equation*}

For $M^{-1/2}$, we use the power series $ (I-A)^{-1/2}=\sum c_{n}A^{n}$.  As $0\le c_{n}\lesssim n^{1/2}$, the argument is the same. 
Finally, if $M$ is not {positive} self-adjoint, then we use $M^{-1}=M^*(MM^*)^{-1}$ and the remark that on a 1-separated set $d\ge d^s$. 
\end{proof}

{
\begin{proof}[Proof of Proposition \ref{prop:matrices} and Theorem~\ref{thm:biortho}]
Lemma~\ref{lemma2} gives us the desired decay for $M_{k}^{-1}$ and $M_{k}^{-1/2}$, as claimed in Proposition~\ref{prop:matrices}, with uniform control of the constant with respect to $k$. The decay and regularity of the $\tilde{s}^k_\alpha$ and $\phi^k_\alpha$, as asserted by Theorem~\ref{thm:biortho}, then follow from the support and regularity of the $s^k_\alpha$ and the decay of the matrix coefficients. We skip details.
\end{proof}
}

%

\begin{remark}\label{rem:ds}
Let  $0<s\le 1$. If $d$ is a quasi-distance with quasi-triangle constant $A_{0}$, then $d^s$ is a quasi-distance with quasi-triangle constant $A_{0}^s$.  If $d$ is geometrically doubling {with constant $N$, then so is $d^s$ with constant at most $N^{\lceil 1/s\rceil}$, where $\lceil\ \rceil$ is the `rounding-up' to the next integer. If $\mu$ is a doubling measure with respect to the balls for $d$ with constant $c_\mu$, then it is with respect to the balls for $d^s$ with constant at most $c_\mu^{\lceil 1/s\rceil}$. Thus both previous lemmata apply to $d^s$ with constants that depend only on the original constants and $s$.} We shall use this remark later.
\end{remark}

\section{Spline wavelets} 

We are now prepared for the construction of an orthonormal basis of $L^2(\mu)$, consisting of wavelets $\psi^k_\alpha$ with similar decay and regularity properties as with the spline systems. We follow an algorithm from Meyer \cite{M}.

Fix $k\in \Z, k\ge k_{0}$, where $k_0$ is as in Theorem~\ref{thm:multires}. Recall the operator
\begin{equation*}
   U_{k+1}^{-1}: f\mapsto \{f(x^{k+1}_{\beta})\sqrt{\mu^{k+1}_{\beta}}\}_{\beta} 
\end{equation*}
is an isomorphism from $V_{k+1}$ onto $\ell^2(\mathscr{X}^{k+1})$. Denote by $Y_{k}$ the inverse image of the subspace of $\ell^2(\mathscr{X}^{k+1})$ sequences vanishing on $\mathscr{X}^{k}$; this subspace is naturally identified with $\ell^2(\mathscr{Y}^k)$, where we recall the notation $\mathscr{Y}^k:=\mathscr{X}^{k+1}\setminus\mathscr{X}^k$. Clearly, 
 $V_{k}\oplus Y_{k}= V_{k+1}$    topologically. Consider the orthogonal (in $L^2(\mu)$) {complement} $W_{k}$ of $V_{k}$ in $V_{k+1}$.  Then the restriction  to $Y_{k}$ of the orthogonal projection $Q_{k}$ onto $W_{k}$  is an isomophism onto $W_{k}$.  Now, identifying with $\mathscr{Y}^k$ the set of indices $\beta$ corresponding to $x^{k+1}_\beta\in\mathscr{X}^{k+1}\setminus\mathscr{X}^k=\mathscr{Y}^k$,
  the collection  $\{s^{k+1}_{\beta}\}_{\beta\in \mathscr{Y}^{k}}$ forms an unconditional basis of $Y_{k}$ and its image under $Q_{k}$ is an unconditional basis of $W_{k}$. A representation of $Q_{k}f$ when $f\in V_{k+1}$ is
\begin{equation}\label{eq:Qkf}
 Q_{k}f= f- \sum_{\alpha\in \mathscr{X}^{k}} \langle f, \tilde s^k_{\alpha}\rangle_{L^2(\mu)} s^k_{\alpha}
 =f- \sum_{\alpha\in \mathscr{X}^{k}} \langle f,  s^k_{\alpha}\rangle_{L^2(\mu)} \tilde{s}^k_{\alpha},
\end{equation}
 because the sum is the orthogonal projection of $f$ onto $V_{k}$.  Hence, the pre-wavelets 
 $$\tilde \psi^k_{\beta}(x):=Q_{k}s^{k+1}_{\beta}(x)$$
 have the $L^\infty$-normalized  exponential decay $C \exp\big( -\gamma ( \delta^{-k}{d(y^{k}_{\beta},x)})^s\,\big)$, where $y^k_{\beta}:=x^{k+1}_{\beta}\in\mathscr{Y}^k$. (Remark that they are normalized as the splines $s^{k+1}_{\beta}$.) Finally, one can orthonormalize them in $L^2(\mu)$ following the procedure of Section \ref{sect:biortho} applied to the positive self-adjoint matrix
 $$
   {\tilde{M}(\alpha,\beta):=}\frac{\langle  \tilde \psi^k_{\alpha},  \tilde \psi^k_{\beta}\rangle}{\sqrt{\mu_{\alpha}^{k+1}{\mu_{\beta}^{k+1}}}}
 $$
 indexed by  $\mathscr{Y}^{k}\times \mathscr{Y}^{k}$; that is, we define
 \begin{equation*}
  \psi^k_\alpha(x):=\sum_{\beta\in\mathscr{Y}^k}\tilde{M}^{-1/2}(\alpha,\beta)\frac{\tilde\psi^k_\beta(x)}{\sqrt{\mu^{k+1}_\beta}}.
\end{equation*}
Note that with the notation $y^k_{\beta}=x^{k+1}_{\beta}$, we have $\mu^{k+1}_{\beta}=\mu(B(y^k_{\beta},\delta^{k+1}))\sim \mu(B(y^k_{\beta},\delta ^k))$. 
The point-set $\mathscr{Y}^k$ is a 1-separated set for $d_{k}^s (x,y):=({\delta ^{-k-1}}{d(x, y)})^s$ and the matrix $\tilde{M}$ has the exponential decay, {as required in Lemma~\ref{lemma2}, with respect to $d_{k}^s$. By Lemma~\ref{lemma2} and Remark~\ref{rem:ds}, $\tilde{M}^{-1/2}$ has decay
\begin{equation*}
  \abs{\tilde{M}^{-1/2}(\alpha,\beta)}
  \lesssim\exp\big(-\gamma' d_k^s(x^k_{\alpha},x^k_{\beta})^{\tilde s}\,\big)=\exp\big(-\gamma(\delta^{-k} d(x^k_{\alpha},x^k_{\beta}))^{s\tilde{s}}\,\big),
\end{equation*}
where $\tilde{s}=(1+\log_2 A_0^s)^{-1}=(1+(1+\log_2 A_0)^{-1}\log_2 A_0)^{-1}=(1+\log_2 A_0)(1+2\log_{2} A_{0})^{-1}$ or $\tilde s=1$ if $d$ is Lipschitz-continuous.

This yields an orthonormal basis $\psi^k_{\beta}(x)$, $y^k_{\beta}\in\mathscr{Y}^k$, of $W_{k}$ having  the $L^2$ normalized   exponential decay 
\begin{equation*}
  |\psi^k_{\beta}(x)|\le \frac{C}{\sqrt{ \mu(B(y^k_{\beta},\delta ^k))}} \exp\big( -\gamma (\delta^{-k}{d(y^{k}_{\beta},x)})^{a}\,\big), \qquad
   a={ s\tilde{s}}.
\end{equation*}
  Gathering the construction for all $k$ (and adding constants if $X$ is bounded), one obtains an orthonormal basis of $L^2(\mu)$ made of spline wavelets. 
 
 \begin{theorem}\label{th:wavelets} Let $(X,d,\mu)$ be any space of homogeneous type with quasi-triangle constant $A_0$, and $a:=(1+2\log_2 A_0)^{-1}$ or $a:=1$ if $d$ is Lipschitz-continuous.
There exists an orthonormal basis $\psi^k_{\beta}$, $k\in \Z$ (and $k\ge k_{0}$ if $X$ is bounded), $y^k_{\beta}\in\mathscr{Y}^k$,  of $L^2(\mu)$  (or the orthogonal space to constants  if $X$ is bounded)  having exponential decay
 \begin{equation*}
  \abs{\psi^k_{\beta}(x)}\leq \frac{C}{\sqrt{ \mu(B(y^k_{\beta},\delta ^k))}} \exp\big( -\gamma (\delta^{-k}{d(y^{k}_{\beta},x)})^{a}\, \big)
,
\end{equation*}
H\"older-regularity
\begin{equation*}
  \abs{\psi^k_{\beta}(x)-\psi^k_{\beta}(y)}\leq \frac{C}{\sqrt{ \mu(B(y^k_{\beta},\delta ^k))}}
    \Big(\frac{d(x,y)}{\delta^k}\Big)^{\eta} \exp\big( -\gamma (\delta^{-k}{d(y^{k}_{\beta},x)})^{a}\, \big)
,\qquad
    d(x,y)\leq\delta^k,
\end{equation*}
 and vanishing mean
 \begin{equation*}
 \int_{X} \psi^k_{\beta}(x)\, d\mu(x)=0, \quad k\in \Z, k \ge k_{0}, y^k_{\beta}\in\mathscr{Y}^k.
\end{equation*}
  \end{theorem}

\begin{proof} It remains to see vanishing mean and  regularity.

Using $\int\tilde{s}^k_\alpha\ud\mu =1$ and $\sum_\alpha s^k_\alpha\equiv 1$, we deduce from \eqref{eq:Qkf} for $f\in L^1(\mu)$ that
\begin{equation*}
  \int Q_k f\ud\mu=\int f\ud\mu-\sum_\alpha\langle f,s^k_\alpha\rangle \int\tilde{s}^k_\alpha\ud\mu
  =\int f\ud\mu-\sum_\alpha\int f s^k_\alpha\ud\mu
  =\int f\ud\mu-\int f\ud\mu=0.
\end{equation*}
Since the pre-wavelets $\tilde\psi^k_\alpha$ lie in the range of $Q_k$ by definition, we have $\int\tilde\psi^k_\alpha\ud\mu=0$, and the same result for the wavelets follows from the convergent series representation of $\psi^k_\alpha$ in terms of the $\tilde\psi^k_\beta$.
   
 As for the regularity, recall the smoothness of the splines: $\abs{s^k_{\alpha}(x)-s^k_{\alpha}(y)}\leq C(d(x,y)/\delta^k)^{\eta}$. This implies for the pre-wavelets the estimate \begin{equation*}
\begin{split}
  &\abs{\tilde\psi^k_{\beta}(x)-\tilde\psi^k_\beta(y)}
  \leq\abs{s^{k+1}_{\beta}(x)-s^{k+1}_{\beta}(y)}
  +\sum_{\alpha}\abs{\pair{s^{k+1}_{\beta}}{\tilde{s}^k_{\alpha}}}\abs{s^k_{\alpha}(x)-s^k_{\alpha}(y)} \\
  &\leq C\Big(\frac{d(x,y)}{\delta^k}\Big)^{\eta}1_{B(y^k_{\beta},c\delta^k)}(x)
  +\sum_{\alpha}C \exp\big (-\gamma ( d(y^k_\beta,x^k_\alpha)/\delta^k)^s\, \big)\Big(\frac{d(x,y)}{\delta^k}\Big)^{\eta}1_{B(x^k_{\alpha},C\delta^k)}(x) \\
  &\leq C\Big(\frac{d(x,y)}{\delta^k}\Big)^{\eta}\exp\big(-\gamma (d(x,y^{k}_\beta)/\delta^k)^s\, \big).
\end{split}
\end{equation*}
For the wavelets, finally, we have
\begin{equation*}
\begin{split}
  \abs{\psi^k_{\alpha}(x)-\psi^k_\alpha(y)}
  &\leq C\sum_\beta\frac{\exp\big({-\gamma (d(y^{k}_{\alpha},y^{k}_{\beta})/\delta^{k+1}})^a\, \big)}{\sqrt{ \mu(B(y^k_{\beta},\delta ^k))}}\abs{\tilde\psi^k_{\beta}(x)-\tilde\psi^k_\beta(y)} \\
 &\leq C\sum_\beta\frac{\exp\big({-\gamma (d(y^{k}_{\alpha},y^{k}_{\beta})/\delta^{k+1}})^a\, \big)}{\sqrt{ \mu(B(y^k_{\beta},\delta ^k))}}
 \Big(\frac{d(x,y)}{\delta^k}\Big)^{\eta}\exp\big(-\gamma (d(x,y^{k}_\beta)/\delta^k)^s\, \big) \\
 &\leq \frac{C}{\sqrt{ \mu(B(y^k_{\alpha},\delta ^k))}}\Big(\frac{d(x,y)}{\delta^k}\Big)^{\eta}\exp\big(-\gamma (d(x,y^{k}_\alpha)/\delta^k)^a\, \big).
\end{split}
\end{equation*}
Note that the value of $\gamma>0$ changes from line to line in these computations. Also we used that $a\le s$  and a variant of Lemma \ref{lemma1} for the quasi-distance $d(x,y)/\delta^k$ on $X$.
\end{proof}

\begin{remark}
The construction and also the next sections  suggest that the label $k+1$ would be more appropriate than $k$   for $W_{k}$, the wavelets $\psi^k_{\alpha}$ and their scale $\delta ^k$, because its keeps closer to the definition of the point sets $\mathscr{Y}^k$, a subset of $\mathscr{X}^{k+1}$, which will take an important role.  We have kept the  wavelet community notation as in \cite{M2}.
\end{remark}
 
\section{Technical estimates related to vanishing annuli}
 
 We break the development of the wavelet theory with this technical section, which will provide us with useful estimates to streamline the subsequent presentation. A basic difficulty related to general spaces of homogeneous type, as opposed to those with the reverse doubling property, is the possible existence of arbitrarily large empty annuli $B(x,R)\setminus B(x,r)=\varnothing$. This leads to a certain dichotomy: locally, we either have the reverse doubling estimate, or the vanishing of a certain annulus, both of which provide certain control, which we need to exploit in different ways. This is quantified in the following. (The next result is well-known, cf. \cite[Remark 1.2]{HMY2}, but we include it here for completeness, since we need it to derive some consequences which appear to be new.)
 
\begin{lemma}\label{lem:emptyAnnulus}
For every $x\in X$ and $R>r>0$, at  least  one of the following alternatives holds:
\begin{equation*}
  \mu(B(x,R))\geq(1+\varepsilon)\mu(B(x,r))\qquad\text{or}\qquad B(x,\frac{1}{2A_{0}}R)\setminus B(x,2A_{0}r)=\varnothing,
\end{equation*}
where $\varepsilon:=1/C_\mu(3A_{0}^2)$, and $C_\mu(t)$ is the smallest constant such that $\mu(tB)\leq C_\mu(t)\mu(B)$ for all balls $B\subseteq X$. 
\end{lemma}

In particular, as $R\to\infty$, we observe the following (again well-known fact): If there is a ball $B(x,r)$ such that $\mu(X)<(1+C_\mu(3A_{0}^2)^{-1})\mu(B(x,r))$, then $X=B(x,2A_0r)$. If $\mu(X)<\infty$, such a ball always exists, and hence $\diam(X)<\infty$.

\begin{proof}
Suppose that the annulus is nonempty, and let $y\in B(x,\frac{1}{2A_{0}}R)\setminus B(x,2A_{0}r)$. Let $\rho:=\frac{1}{2A_{0}}d(x,y)\geq r$. We claim that
\begin{equation*}
  B(y,\rho)\subseteq B(x,R)\setminus B(x,r).
\end{equation*}
Indeed, if $z\in B(y,\rho)$, then
\begin{equation*}
  d(z,x)\leq A_{0}d(z,y)+A_{0}d(y,x)<A_0\rho+A_{0}d(x,y)\leq R,
\end{equation*}
while
\begin{equation*}
  d(z,x)\geq\frac{1}{A_{0}}d(y,x)-d(y,z)>\frac{1}{A_{0}}d(x,y)-\rho=\frac{1}{2A_{0}}d(x,y)\geq r.
\end{equation*}

We also claim that
\begin{equation*}
  B(x,r)\subset B(y,3A_{0}^2 \rho).
\end{equation*}
Indeed, if $w\in B(x,r)$, then
\begin{equation*}
  d(w,y)\leq A_{0}d(w,x)+A_{0}d(x,y)<A_{0}r+2A_{0}^2\rho\leq 3A_{0}^2\rho.
\end{equation*}

Now it follows that
\begin{equation*}
  \mu(B(x,R))-\mu(B(x,r))
  \geq\mu(B(y,\rho))
  \geq\frac{1}{C_\mu(3A_{0}^2)}\mu(B(y,3A_{0}^2\rho))
  \geq\frac{1}{C_\mu(3A_{0}^2)}\mu(B(x,r)).\qedhere
\end{equation*}
\end{proof}

The following lemma relates the mentioned dichotomy to the distribution of the dyadic point sets $\mathscr{Y}^k$:
  
\begin{lemma}\label{lem:kjSeq}
For every $x\in X$ and $r>0$, there exists a decreasing sequence, finite or infinite, of integers $\{k_j\}_{j=0}^J$ such that $r\leq\delta^{k_0}<\delta^{k_1}<\ldots$ such that
\begin{equation*}
  V(x,\delta^k)\gtrsim(1+\varepsilon)^j V(x,r)\quad\text{and}\quad
  d(x,\mathscr{Y}^k)+\delta^k\gtrsim\delta^{k_{j+1}}\quad\text{if }k_j\geq k>k_{j+1},
\end{equation*}
where we intepret $k_{J+1}:=-\infty$ if $J<\infty$.
\end{lemma}

\begin{proof}
Let $k(0)$ be the largest integer with $\delta^{k(0)}\geq r$, and let $k(j+1)$ be the largest integer with $V(x,\delta^{k(j+1)})\geq(1+\varepsilon)V(x,\delta^{k(j)})$. (Note that the sequence terminates if and only if $\mu(X)<\infty$.) Thus $V(x,\delta^{k(j+1)+1})<(1+\varepsilon)V(x,\delta^{k(j)})$, and hence
\begin{equation*}
  B(x,\frac{1}{2A_{0}}\delta^{k(j+1)+1})=B(x,2A_{0}\delta^{k(j)}).
\end{equation*}
For $k\geq k(j+1)+2$, the ball on the left contains at least one element of $\mathscr{X}^k$. For $k\leq k(j)-1$, the ball on the right contains at most one element of $\mathscr{X}^k$. Since the balls are equal, for $k(j+1)+2\leq k\leq k(j)-1$, the ball contains exactly one element of $\mathscr{X}^k$. So if $k$ and $k+1$ are both in this range, i.e., if $k(j+1)+2\leq k\leq k(j)-2$, then the intersections of the ball with $\mathscr{X}^k$ and $\mathscr{X}^{k+1}$ coincide; hence there is no point of $\mathscr{Y}^k=\mathscr{X}^{k+1}\setminus\mathscr{X}^k$ in the ball. Thus $d(x,\mathscr{Y}^k)\geq \frac{1}{2A_{0}}\delta\cdot\delta^{k(j+1)}$. On the other hand, if $k\in\{k(j+1)-1,k(j+1),k(j+1)+1\}$, then clearly $\delta^k\geq\delta\cdot\delta^{k(j+1)}$, and hence
\begin{equation*}
  d(x,\mathscr{Y}^k)+\delta^k\geq\frac{\delta}{2A_{0}}\delta^{k(j+1)}\qquad\text{if } k(j+1)-1\leq k\leq k(j)-2.
\end{equation*}
Also, for $k$ in the same range, we have
\begin{equation*}
  V(x,\delta^k)
  \geq V(x,\delta^{k(j)})\geq (1+\varepsilon)^j V(x,\delta^{k(0)})\geq (1+\varepsilon)^j V(x,r).
\end{equation*}
So the claim follows by relabeling $k_{j+1}:=k(j)-2$ for $j\in\N$ and $k_0:=k(0)$. 
\end{proof}

Sums of the following type appear in connection with the wavelets:

\begin{lemma}\label{lem:sumLargeBalls} 
For all $x\in X$ and $r,\nu, a>0$, we have
\begin{equation*}
  \sum_{k:\delta^k\geq r}V(x,\delta^k)^{-\nu}\exp\big(-\gamma(\delta ^{-k}{d(x,\mathscr{Y}^k)})^a\,  \big)\lesssim V(x,r)^{-\nu}.
\end{equation*}
\end{lemma}

\begin{proof}
Let $k_j$ be the sequence as provided by Lemma~\ref{lem:kjSeq}. Then
\begin{equation*}
\begin{split}
 \sum_{k:\delta^k\geq r} &V(x,\delta^k)^{-\nu}\exp\big(-\gamma(\delta ^{-k}{d(x,\mathscr{Y}^k)})^a\,  \big)\\
 &=\sum_{j=0}^J\sum_{k:k_j\geq k>k_{j+1}}
    V(x,\delta^k)^{-\nu}\exp\big(-\gamma(\delta ^{-k}{d(x,\mathscr{Y}^k)})^a\,  \big) \\
 &\lesssim\sum_{j=0}^J\sum_{k:k_j\geq k>k_{j+1}}(1+\varepsilon)^{-j\nu}V(x,r)^{-\nu} 
    \exp\big(-\gamma\delta^{a(k_{j+1}-k)}\, \big) \\
 &\leq V(x,r)^{-\nu} \Big(\sum_{j=0}^J (1+\varepsilon)^{-j\nu}\Big)
     \Big(\sum_{m=1}^{\infty} \exp(-\gamma\delta^{-ma})\Big) \lesssim V(x,r)^{-\nu}.\qedhere
\end{split}
\end{equation*}
\end{proof}

\section{Technical estimates involving the wavelets}\label{sec:techWavelets}\label{sec:technical}

Actually, the estimates here are valid for any family of functions $\psi^k_\alpha$ which satisfy the same size and regularity estimates as the wavelets. This includes the condition that $\psi^k_\alpha$ be concentrated around the point $y^k_\alpha\in\mathscr{Y}^k$, and the structure of the point sets $\mathscr{Y}^k$ is important for some of the following estimates.

\begin{lemma}\label{lem:sumAlpha} Let $a$ as in Theorem \ref{th:wavelets}.
For a fixed $k\in\Z$,
\begin{equation*}
  \sum_{\alpha\in\Lambda_k}\abs{\psi^k_\alpha(x)\psi^k_\alpha(y)}
  \leq\frac{C}{V(x,\delta^k)}\exp\big(-\gamma (\delta ^{-k}{d(x,\mathscr{Y}^k)})^a \, \big)\exp\big(-\gamma (\delta ^{-k}{d(x,y)})^a \, \big),
\end{equation*}
and, for $d(x,x')<\frac{1}{2A_0}d(x,y)$,
\begin{equation*}
\begin{split}
    \sum_{\alpha\in\Lambda_k} &\abs{[\psi^k_\alpha(x)-\psi^k_\alpha(x')]\psi^k_\alpha(y)} \\
  &\leq\frac{C}{V(x,\delta^k)}\min\Big\{1,\Big(\frac{d(x,x')}{\delta^k}\Big)^{\eta}\Big\}\exp\big(-\gamma (\delta ^{-k}{d(x,\mathscr{Y}^k)})^a \, \big)\exp\big(-\gamma (\delta ^{-k}{d(x,y)})^a \, \big).
\end{split}
\end{equation*}
\end{lemma}
 
Note that if, instead, we have a family of functions $\varphi^k_\alpha$ corresponding to the points $x^k_\alpha\in\mathscr{X}^k$ instead of $y^k_\alpha\in\mathscr{Y}^k$, we get exactly the same estimate but with $d(x,\mathscr{X}^k)$ in place of $d(x,\mathscr{Y}^k)$ in the result. Since  $d(x,\mathscr{X}^k)\leq 2A_{0}\delta^k$, the first exponential factor is roughly $1$, and may be dropped.

\begin{proof}
By the doubling condition and the quasi-triangle inequality, we estimate
\begin{equation*}
\begin{split}
   \abs{\psi^k_\alpha(x)\psi^k_\alpha(y)}
   &\leq \frac{C}{\sqrt{\mu^{k+1}_\alpha}}\exp\big(-\gamma (\delta ^{-k}{d(x,y^k_{\alpha})})^a \, \big)
   \frac{C}{\sqrt{\mu^{k+1}_\alpha}}\exp\big(-\gamma (\delta ^{-k}{d(y,y^k_\alpha)})^a \, \big) \\
   &\leq \frac{C}{V(x,\delta^k)}\exp\big(-\gamma (\delta ^{-k}{d(x,y^k_{\alpha})})^a \, \big)\exp\big(-\gamma (\delta ^{-k}{d(x,y)})^a \, \big),
\end{split}
\end{equation*}
and the sum over $\alpha\in \mathscr{Y}^k$ of the second factor is dominated by $\exp\big(-\gamma (\delta ^{-k}d(x,\mathscr{Y}^k))^a \, \big)$.

If $\delta^k\geq d(x,x')$, then
\begin{equation*}
\begin{split}
  \abs{[\psi^k_\alpha(x)-\psi^k_\alpha(x')]\psi^k_\alpha(y)}
  &\leq\frac{C}{\mu^{k+1}_\alpha}\Big(\frac{d(x,x')}{\delta^k}\Big)^{\eta}\exp\big(-\gamma (\delta ^{-k}{d(x,y^k_{\alpha})})^a \, \big)
    \exp\big(-\gamma (\delta ^{-k}{d(y,y^k_{\alpha})})^a \, \big) \\
  &\leq\frac{C}{V(x,\delta^k)}\Big(\frac{d(x,x')}{\delta^k}\Big)^{\eta}\exp\big(-\gamma (\delta ^{-k}{d(x,y^k_{\alpha})})^a \, \big)\exp\big(-\gamma (\delta ^{-k}{d(x,y)})^a \, \big),
\end{split}
\end{equation*}
and we may similarly sum over $\alpha\in\mathscr{Y}^k$. For $\delta^k<d(x,x')\leq d(x,y)$, we just use the {quasi-triangle} inequality and the first estimate of the lemma to both terms.
\end{proof}
 
\begin{lemma}\label{lem:sumAlphaK}
\begin{equation*}
  \sum_{k,\alpha}\abs{\psi^k_\alpha(x)\psi^k_\alpha(y)}
  \leq\frac{C}{V(x,y)}.
\end{equation*}
\end{lemma}
 
\begin{proof}
By Lemma~\ref{lem:sumAlpha}, we have
\begin{equation*}
\begin{split}
   \sum_{k,\alpha}\abs{\psi^k_\alpha(x)\psi^k_\alpha(y)}
   &\leq\sum_{k:\delta^k\geq d(x,y)}\frac{C}{V(x,\delta^k)}\exp\big(-\gamma (\delta ^{-k}{d(x,\mathscr{Y}^k)})^a \, \big) \\
   &\qquad+\sum_{k:\delta^k<d(x,y)}\frac{C}{V(x,\delta^k)}\exp\big(-\gamma (\delta ^{-k}{d(x,y)})^a \, \big),
\end{split}
\end{equation*}
where the first part has the correct bound by Lemma~\ref{lem:sumLargeBalls}. For the second part, we have
\begin{equation*}
\begin{split}
  \sum_{k:\delta^k<d(x,y)} &\frac{C}{V(x,\delta^k)}\exp\big(-\gamma (\delta ^{-k}{d(x,y)})^a \, \big) \\
  &\leq\sum_{k:\delta^k<d(x,y)}\frac{C}{V(x,y)}\Big(\frac{d(x,y)}{\delta^k}\Big)^M\exp\big(-\gamma (\delta ^{-k}{d(x,y)})^a \, \big) \\
  &\leq\frac{C}{V(x,y)}\sum_{m=0}^{\infty}\delta^{-mM}\exp(-\gamma\delta^{-ma})\leq\frac{C}{V(x,y)}.\qedhere
\end{split}
\end{equation*}
\end{proof}
 
 \begin{lemma}\label{lem:sumDifference}
\begin{equation*}
  \sum_{k,\alpha}\abs{[\psi^k_\alpha(x)-\psi^k_\alpha(x')]\psi^k_\alpha(y)}
  \leq\frac{C}{V(x,y)}\Big(\frac{d(x,x')}{d(x,y)}\Big)^{\eta},\qquad
  d(x,x')<\frac{1}{2A_0}d(x,y).
\end{equation*}
\end{lemma}
 
\begin{proof}
By the second estimate of Lemma~\ref{lem:sumAlpha},
\begin{equation*}
\begin{split}
  \sum_{k,\alpha} &\abs{[\psi^k_\alpha(x)-\psi^k_\alpha(x')]\psi^k_\alpha(y)} \\
  &\leq\sum_{k:\delta^k\geq d(x,y)}\frac{C}{V(x,\delta^k)}\Big(\frac{d(x,x')}{d(x,y)}\Big)^{\eta}\exp\big(-\gamma (\delta ^{-k}{d(x,\mathscr{Y}^k)})^a \, \big)  \\
  &\qquad+\sum_{k:d(x,x')\leq \delta^k\leq d(x,y)}\frac{C}{V(x,y)}
   \Big(\frac{d(x,y)}{\delta^k}\Big)^{M+\eta}\Big(\frac{d(x,x')}{d(x,y)}\Big)^{\eta}\exp\big(-\gamma (\delta ^{-k}{d(x,y)})^a \, \big) \\
  &\qquad+\sum_{k:\delta^k\leq d(x,x')}\frac{C}{V(x,y)}
   \Big(\frac{d(x,y)}{\delta^k}\Big)^{M}\exp\big(-\gamma (\delta ^{-k}{d(x,y)})^a \, \big).
\end{split}
\end{equation*}
The first two parts contain the factor $(d(x,x')/d(x,y))^{\eta}$, and the rest is bounded by $C/V(x,y)$ according to Lemma~\ref{lem:sumLargeBalls} and $\sum_{m=0}^{\infty}\delta^{-m(M+\eta)}\exp(-\gamma\delta^{-ma})\leq C$. For the last term we even obtain the bound $C/V(x,y)\cdot\exp\big(-\gamma (d(x,y)/d(x,x'))^a\, \big)\leq C/V(x,y)\cdot (d(x,x')/d(x,y))^K$ for any $K$.
\end{proof}

 \section{Littlewood--Paley decomposition and $L^p$ theory}

Recall that $Q_{k}$ is the orthogonal projector onto $W_{k}$ and set also $P_{k}$ the orthogonal projector onto $V_{k}$. These operators will provide us with a new regular Littlewood--Paley decomposition for spaces of homogeneous type. The following lemma describes the kernels of these operators:

\begin{lemma} The kernel $P_{k}(x,y)$ of $P_{k}$ is  symmetric in $x,y$ and  has estimates
$$|P_{k}(x,y)|\le  \frac{C}{\sqrt{\mu(B(x,\delta^k))\mu(B(y,\delta ^k))}} \exp\big(-\gamma (\delta ^{-k}{d(x,y)})^s \, \big)
$$
$$|P_{k}(x,y)-P_{k}(x,y')|\le C \left({\frac{d(y,y')}{\delta ^{k}}}\right)^\eta
\left ( \frac{\exp\big(-\gamma (\delta ^{-k}{d(x,y)})^s \, \big)}{\sqrt{\mu(B(x,\delta^k))\mu(B(y,\delta ^k))}}+ \frac{\exp\big(-\gamma (\delta ^{-k}{d(x,y')})^s \, \big)}{\sqrt{\mu(B(x,\delta^k))\mu(B(y',\delta ^k))}}\right)$$
for some $C,\gamma$ and all $x,y,y'\in X$ and $k\in \Z$ (with $k\ge k_{0}$ if $X$ is bounded). 
Moreover
$$
\int_{X}P_{k}(x,y)d\mu(x)=1.
$$
The kernel $Q_{k}(x,y)$ of $Q_{k}$ is  symmetric in $x,y$ and  has similar estimates with {$s$ changed to $a$}, the additional exponential factor
\begin{equation*}
  \exp\big(-\gamma(\delta^{-k}d(x,\mathscr{Y}^k))^a \, \big),
\end{equation*}
and the cancellation condition
$$
\int_{X}Q_{k}(x,y)d\mu(x)=0.
$$
\end{lemma}

Note that there are many ways to express  the denominator factors in the kernels up to changing the constants $C,\gamma$ in the numerator factors especially thanks to the exponential decay. In particular, one possible expression shows that the system of operators $P_{k}$ is an `Approximation of The Identity' in the sense of  \cite[Definition 2.2]{HMY1} (this part of that paper does not use the Reverse Doubling property), and the properties listed there hold.

 \begin{theorem} The spline-wavelet representation  yields a  decomposition of Littlewood-Paley type with H\"older-continuous kernels. 

\end{theorem}

\begin{proof}
If $X$ is unbounded and $f\in L^2(X)$ the converging series
 $$
 f= \sum_{k=-\infty}^\infty Q_{k}f= \sum_{k=-\infty}^\infty Q_{k}^2f
 $$
 is an homogeneous Littlewood-Paley decomposition.  Also one can truncate at any level  since $P_{k+1}=P_{k}+Q_{k}$ and write for any $\ell$
 $$
 f= P_{\ell}^2f+ \sum_{k=\ell}^\infty Q_{k}^2f
 $$
 which is  an inhomogeneous Littlewood-Paley decomposition. 
 This decomposition is the one used if $X$ is bounded with $\ell=k_{0}$. In that case, the term $P_{\ell}^2f$ is a constant. 
 \end{proof} 
 
Observe that the sum $ \sum_{k} Q_{k}f$ can also be rewritten as a ``discrete'' Littlewood-Paley decomposition $\sum_{k,\alpha} \langle f, \psi^k_{\alpha}\rangle \psi^k_{\alpha}$. 
 From there, one can look at convergence for $f$ in various topological spaces, and develop the theory of function spaces. We restrict ourselves to $L^p$ spaces (and BMO in the next section) and leave further developments to the interested reader.

The estimates in Section~\ref{sec:techWavelets} immediately give the following result:

\begin{proposition}\label{prop:waveletsToCZO}
Let $c^k_\alpha$ be arbitrary complex coefficients bounded in absolute value by one. Then the series
\begin{equation*}
  K(x,y)=\sum_{k\in\Z}\sum_{\alpha\in\Lambda_k} c^k_\alpha\psi^k_\alpha(x)\psi^k_\alpha(y)
\end{equation*}
converges absolutely for $x\neq y$ and satisfies
\begin{equation*}
  \abs{K(x,y)}\leq\frac{C}{V(x,y)},\qquad
  \abs{K(x,y)-K(x',y)}\leq\frac{C}{V(x,y)}\Big(\frac{d(x,x')}{d(x,y)}\Big)^{\eta},\quad d(x,x')\leq\frac{1}{2A_{0}}d(x,y),
\end{equation*}
with a similar H\"older-regularity estimate in the second variable.
\end{proposition}

\begin{corollary}
The spline wavelets form an unconditional basis of $L^p(\mu)$ spaces when $1<p<\infty$. 
\end{corollary}

\begin{proof} 
Completeness follows using the convergence properties of the $P_{k}$ in $L^p(\mu)$. It remains to show that  operators $T_{c}$  given $T_{c}(\psi^k_{\alpha})=c^k_{\alpha}\psi^k_{\alpha}$ are uniformly bounded on all $L^p(\mu)$ spaces whenever $c=(c^{k}_{\alpha})$ is a sequence of complex numbers in the unit ball of $\mathbb{C}$. These operators are contractions in $L^2(\mu)$.  Using the regularity of their kernels  as proved in Proposition \ref{prop:waveletsToCZO} and the Calder\'on--Zygmund theorem, $T_{c}$ has weak type (1,1), with uniform bound with respect to $c$. Since $T_{c}^*=T_{\bar c}$, the same applies to $T_{c}^*$. We conclude by interpolation.  
\end{proof}

\section{BMO theory}

Recall that the space $\BMO(\mu)$ of functions of \textbf{bounded mean oscillation} consists of those $b\in L^1_{\loc}(\mu)$ with
\begin{equation*}
  \Norm{b}{\BMO(\mu)}:=\sup_B\inf_c\frac{1}{\mu(B)}\int_B\abs{b-c}\ud\mu<\infty,
\end{equation*}
where the infimum is almost realized by $c=b_B:=\mu(B)^{-1}\int_B b\ud\mu$. { We do not incorporate {the norm of constants as a part of the $\BMO$ norm} even if $X$ is bounded, so that our $\BMO$ is a Banach space modulo constants in both bounded and unbounded cases.} These averages satisfy the following useful estimate:

\begin{lemma}\label{lem:bBs}
Let $B_i=B(x_i,r_i)$, $i=1,2$, be two balls in $X$. Then
\begin{equation*}
  \abs{b_{B_1}-b_{B_2}}\lesssim\Norm{b}{\BMO(\mu)}\Big(1+\log\frac{r_1+r_2+d(x_1,x_2)}{\min\{r_1,r_2\}}\Big).
\end{equation*}
\end{lemma}

\begin{proof}
Suppose first that $r_1\leq r_2$ and $B_1\subset B_2$. Let $B^0:=B_1$ and $B^{i+1}$ be the smallest $2^j B^i$ such that $\mu(B^{i+1})>C_\mu\mu(B^i)$; hence also $\mu(B^{i+1})\leq C_\mu^2\mu(B^i)$. It is easy to check that $\abs{b_{B^i}-b_{B^{i+1}}}\lesssim\Norm{b}{\BMO(\mu)}$. Let $B^k$ be the first $B^i$ with radius bigger than that of $B_2$. Then we also have $\abs{b_{B^k}-b_{B_2}}\lesssim\Norm{b}{\BMO(\mu)}$, and $k\lesssim 1+\log r_2/r_1$.

In the general case, we may choose an auxiliary ball $B_3$ of radius $r_3\sim r_1+r_2+d(x_1,x_2)$ which contains both $B_1$ and $B_2$, and apply the earlier consideration to
\begin{equation*}
  \abs{b_{B_1}-b_{B_2}}\leq\abs{b_{B_1}-b_{B_3}}+\abs{b_{B_3}-b_{B_2}}.\qedhere
\end{equation*}
\end{proof}

\begin{corollary}
Let $f$ be a function with exponential decay $\abs{f(x)}\leq C\exp(-c\cdot d(x,x_0)^a)$ with $C,c,a>0$, and $b\in\BMO(\mu)$. Then the product $f\cdot b$ is integrable over $X$.
\end{corollary}

\begin{proof}
Let $B_n:=B(x_0,n)$. Clearly $f$, and hence $f\cdot b_{B_1}$ is integrable, so we consider $f\cdot(b-b_{B_1})$. By Lemma~\ref{lem:bBs},
\begin{equation*}
  \int_{B_n}\abs{b-b_{B_1}}\ud\mu
  \leq\int_{B_n}\abs{b-b_{B_n}}\ud\mu+\mu(B_n)\abs{b_{B_n}-b_{B_1}}
  \lesssim\mu(B_n)(1+\log n)\Norm{b}{\BMO(\mu)},
\end{equation*}
and $\mu(B_n)\lesssim n^M\mu(B_0)$ by doubling. Hence $\Norm{1_{B_n}(b-b_{B_0})}{L^1(\mu)}$ grows at most polynomially in $n$, whereas $\Norm{1_{B_n\setminus B_{n-1}}f}{L^{\infty}(\mu)}$ decays exponentially in $n$. Thus $f\cdot(b-b_{B_0})\in L^1(\mu)$.
\end{proof}

This implies in particular that the wavelet coefficients $(b,\psi^k_\alpha)=\int_X b\cdot\psi^k_\alpha\ud\mu$ are well-defined for $b\in\BMO(\mu)$. The following injectivity property is somewhat technical, and we postpone its proof to Appendix~\ref{app:BMO}:

\begin{proposition}\label{prop:uniquenessBMO}
Suppose that $b\in\BMO(\mu)$ satisfies $(b,\psi^k_\alpha)=0$ for all $k\in\Z$ and $\alpha\in\mathscr{Y}^k$. Then $b$ is equal to a constant.
\end{proposition}

If $X$ is bounded, one should mention that $k\ge k_{0}$, but in fact, notice that $\mathscr{Y}^k=\emptyset$ and so $\psi^k_{\alpha}$ does not even exist  when $k<k_{0}$. We thus do not need to distinguish further between $X$ bounded or not.

Note that some nontrivial a priori size condition on $b$ is in general necessary for such a conclusion. For example, if the $\psi^k_\alpha$ are regular wavelets on $\R^d$, then $(\psi^k_\alpha,P)=0$ for all polynomials $P$ of degree lower than the regularity of the wavelets.
  
We say that a sequence $\{b^k_\alpha\}_{k\in\Z,\alpha\in\mathscr{Y}^k}$ is a \textbf{Carleson sequence} if
\begin{equation*}
  \Norm{\{b^k_\alpha\}_{k,\alpha}}{\operatorname{Car}}
  :=\sup_{\ell\in\Z,\beta\in\mathscr{X}^{\ell}}
  \Big(\frac{1}{\mu(Q^\ell_\beta)}\sum_{\substack{k\in\Z,\alpha\in\mathscr{Y}^k\\(k+1,\alpha)\leq(\ell,\beta)}}\abs{b^k_\alpha}^2\Big)^{1/2}<\infty.
\end{equation*}
Pay attention to the fact that the supremum runs over all $\ell$ and $\beta\in\mathscr{X}^{\ell}$, which index the dyadic cubes $Q^{\ell}_{\beta}$, whereas the sum runs over $k$ and $\alpha\in\mathscr{Y}^k=\mathscr{X}^{k+1}\setminus\mathscr{X}^k$. Via the wavelet decomposition, we obtain an isomorphism between BMO functions and Carleson sequences:

\begin{theorem}\label{th:BMO}
The spaces $\BMO(\mu)/\C$ ($\BMO$ functions modulo constants) and $\operatorname{Car}$ are isomorphic. This isomorphism is realized via $b\mapsto\{(b,\psi^k_\alpha)\}_{k,\alpha}$, with inverse given by
\begin{equation}\label{eq:expandBMO}
  \{b^k_\alpha\}_{k,\alpha}\mapsto
   \sum_{k,\alpha}b^k_\alpha\Big(\psi^k_\alpha-1_{\{j:\delta^j>r_0\}}(k)\psi^k_\alpha(x_0)\Big),
\end{equation}
where the series converges in $L^2_{\loc}(\mu)$ for every $x_0\in X$ and $r_0>0$, and the choices of $x_0$ and $r_0$ only alter the result by an additive constant.
\end{theorem}

The result and its proof are reasonably classical in spirit, but the lack of bounded support of the wavelets somewhat complicates the matters. While exponential decay is intuitively almost as good, one needs to go through certain technicalities if one wants to be careful with convergence issues. We indicate the argument. 

Note that we have formulated the recovery of the $\BMO$ function from its wavelet coefficients using the ``infra-red'' renormalisation of the wavelet series (since the modification on the series appears for large scales which corresponds to small frequencies in the classical Euclidean situation) rather than 
the $H^1$--$\BMO$ duality; this is in contrast to the statement in Meyer's book \cite{M2}, for instance, but this other method would be  certainly doable here as well.

\begin{proof}
Let us first see that $b\mapsto\{(b,\psi^k_\alpha)\}_{k,\alpha}$ maps $\BMO(\mu)/\C$ into Carleson sequences. The injectivity of this mapping is the content of Proposition~\ref{prop:uniquenessBMO}.
Given $(k,\alpha)$, we let $\tilde{B}=B(x^k_\alpha,C\delta^k)$, with some large $C$, be a ball such that $d(Q^k_\alpha,(B^k_\alpha)^c)\gtrsim\delta^k$. We write
\begin{equation*}
   b=(b-b_{\tilde{B}})1_{\tilde{B}}+(b-b_{\tilde{B}})1_{\tilde{B}^c}+b_{\tilde{B}}=:b^1+b^2+b^3.
\end{equation*}
Then
\begin{equation*}
  \sum|(b^1,\psi^\ell_\beta)|^2\leq\|b^1\|_2^2\lesssim\|b\|_{\BMO(\mu)}\mu(Q^k_\alpha)
\end{equation*}
by orthogonality, John--Nirenberg inequality and doubling; $(b^2,\psi^\ell_\beta)$ may be estimated by the decay of the wavelets, and $(b^3,\psi^k_\alpha)=0$.

\subsubsection*{From Carleson sequences to $\BMO$}
Conversely, assume the Carleson condition. Given a ball $B_1=B(x_1,r_1)$, we may rearrange $1_{B_1}$ times the right side of \eqref{eq:expandBMO} as
\begin{equation*}
\begin{split}
  1_{B_1}\sum_{k:\delta^k\leq r_1} &\sum_{\alpha:y^k_\alpha\in CB_1} b^k_\alpha\psi^k_\alpha
  +1_{B_1}\sum_{k:\delta^k\leq r_1}\sum_{\alpha:y^k_\alpha\in (CB_1)^c} b^k_\alpha\psi^k_\alpha
  +1_{B_1}\sum_{k:\delta^k>r_1}\sum_{\alpha}b^k_\alpha(\psi^k_\alpha-\psi^k_\alpha(x_1)) \\
  &-1_{B_1} \sum_{k:r_0<\delta^k\leq r_1}\sum_\alpha b^k_\alpha\psi^k_\alpha(x_0)
    +1_{B_1} \sum_{k:r_1<\delta^k\leq r_0}\sum_\alpha b^k_\alpha\psi^k_\alpha(x_1) \\
   &\qquad+1_{B_1} \sum_{k:\delta^k>r_0\vee r_1}b^k_\alpha(\psi^k_\alpha(x_1)-\psi^k_\alpha(x_0)).
\end{split}
\end{equation*}
Here the last three terms converge uniformly to $1_{B_1}$ times a constant, the first term converges in $L^2(\mu)$ with norm bounded by $\mu(B_1)^{1/2}$ (use the Carleson condition after covering $CB_1$ by boundedly many dyadic cubes of sidelength $\sim r_1$), and the second and third terms can be estimated uniformly by the decay and regularity of the wavelets $\psi^k_\alpha$. This proves the $L^2(\mu)$ convergence on $B_1$, and also the $\BMO$ estimate
\begin{equation*}
  \Norm{1_{B_1}(\tilde{b}-c_{B_1})}{2}\lesssim\mu(B_1)^{1/2},
\end{equation*}
where $\tilde{b}$ stands for the function on the right of \eqref{eq:expandBMO}, and $c_{B_1}$ is the constant produced by the last three terms in the above expansion. Effectively, the same argument also shows the possibility of replacing $(x_0,r_0)$ in \eqref{eq:expandBMO} by $(x_1,r_1)$, only changing the result by a constant.

\subsubsection*{Verifying that the two mappings are inverses to each other}
This is the most technical part of the argument. For a Carleson sequence $\{b^k_\alpha\}_{k,\alpha}$, let $\tilde{b}$ denote the $\BMO(\mu)$ function (as shown in the previous part) on the right of \eqref{eq:expandBMO}. We claim that $(\tilde{b},\psi^k_\alpha)=b^k_\alpha$, which completes the proof. Indeed, this also shows that if $\{b^k_\alpha\}_{k,\alpha}$ arose as $b^k_\alpha=(b,\psi^k_\alpha)$ from some $b\in\BMO(\mu)$, then the function $\tilde{b}$ produced by \eqref{eq:expandBMO} satisfied $(\tilde{b}-b,\psi^k_\alpha)=0$ and hence $\tilde{b}=b+\operatorname{constant}$.

While the claim is formally obvious, due to the orthogonality
\begin{equation*}
 (\psi^\ell_\beta-1_{\{j:\delta^j>r_0\}}({\ell})\psi^\ell_\beta(x_0),\psi^k_\alpha)=\delta_{k\ell}\delta_{\alpha\beta},
\end{equation*}
we need to justify exchanging the order of summation and integration. To this end, let us reinvestigate the convergence of \eqref{eq:expandBMO} where, since $\psi^k_\alpha$ is orthogonal to constants, we may assume that $(x_0,r_0)=(y^k_\alpha,\delta^k)$. Then
\begin{equation*}
  \tilde{b}=\sum_{\theta}\sum_{(\ell+1,\beta)\leq(k{+1},\theta)}b^{\ell}_\beta\psi^\ell_\beta
  +\sum_{\ell:\delta^{\ell}>\delta^k}\sum_\beta b^{\ell}_\beta(\psi^\ell_\beta-\psi^\ell_\beta(y^k_\alpha))=:\tilde{b}^1+\tilde{b}^2.
\end{equation*}
The second part converges pointwise absolutely to a limit of size $\lesssim 1+\log_+ (d(x,y^k_\alpha)/\delta^k)$. Since $\psi^k_\alpha$ decays exponentially, this allows to exchange the summation and integration when computing $(\tilde{b}^2,\psi^k_\alpha)$. Since $\psi^k_\alpha$ is orthogonal to $\psi^{\ell}_\beta$ for $\delta^{\ell}>\delta^k$, as well as to constants, we get $(\tilde{b}^2,\psi^k_\alpha)=0$.

For the first part, we have
\begin{equation*}
  \tilde{b}^1=:\sum_\theta\tilde{b}_\theta=\sum_\theta\sum_\zeta 1_{Q^{k+1}_\zeta}\tilde{b}_\theta,\qquad
  \Norm{1_{Q^{k+1}_\zeta}\tilde{b}_\theta}{2}\lesssim\sqrt{\mu(Q^{k+1}_\zeta)}\exp\big(-\gamma(\delta ^{-k}{d(x^{k+1}_\zeta,x^{k+1}_\theta)})^a\, \big),
\end{equation*}
whereas
\begin{equation*}
  \Norm{1_{Q^{k+1}_\zeta}\psi^k_\alpha}{2}\lesssim\exp\big(-\gamma(\delta ^{-k}{d(x^{k+1}_\zeta,x^{k+1}_\alpha)})^a\, \big).
\end{equation*}
Hence
\begin{align*}
 \sum_\theta &\|\tilde b_\theta\psi^k_\alpha\|_1
  =\sum_\theta\sum_\zeta\|1_{Q^{k+1}_\zeta}\tilde b_\theta\psi^k_\alpha\|_1\\
  &\lesssim\sum_\zeta\sum_\theta\sqrt{\mu(Q^{k+1}_\zeta)}
    \exp\big(-\gamma(\delta ^{-k}{d(x^{k+1}_\zeta,x^{k+1}_\theta)})^a\, \big)\exp\big(-\gamma(\delta ^{-k}{d(x^{k+1}_\zeta,x^{k+1}_\alpha)})^a\, \big) \\
  &\lesssim\sum_\zeta\sqrt{\mu(Q^{k+1}_\zeta)}\exp\big(-\gamma(\delta ^{-k}{d(x^{k+1}_\zeta,x^{k+1}_\alpha)})^a\, \big)\lesssim\sqrt{\mu^{k+1}_\alpha}.
\end{align*}
This allows the first exchange in
\begin{equation*}
  (\tilde{b}^1,\psi^k_\alpha)
  =\sum_\theta(\tilde{b}_\theta,\psi^k_\alpha)
  =\sum_\theta\sum_{(\ell{+1},\beta)\leq(k{+1},\theta)}b^\ell_\beta(\psi^\ell_\beta,\psi^k_\alpha)
  =\sum_\theta\sum_{(\ell{+1},\beta)\leq(k {+1},\theta)}b^\ell_\beta\delta_{\ell k}\delta_{\alpha\beta}
  =b^k_\alpha,
\end{equation*}
where the second exchange follows from the convergence of $\sum_{(\ell{+1},\beta)\leq(k+1,\theta)}b^\ell_\beta\psi^\ell_\beta$ in $L^2(\mu)$, and from $\psi^k_\alpha\in L^2(\mu)$. Altogether, we get $(\tilde{b},\psi^k_\alpha)=b^k_\alpha$, as claimed.
\end{proof}

\section{The $T(1)$ theorem}

{
To illustrate the power of the spline wavelets, we use them to sketch a proof of the $T(1)$ theorem in any space of homogeneous type. Such a result is surely part of the folklore, but surprisingly difficult to find spelled out in complete generality:  the seminal paper of David--Journ\'e--Semmes \cite{DJS}  makes several assumptions on the space, like the small boundary property of balls, and many recent references treat other special cases like Ahlfors--David \cite{DH} or reverse doubling spaces \cite{HMY1}.

The technology to prove the $T(1)$ theorem in a general space of homogeneous type has certainly existed since the work of M.~Christ \cite{Christ}.
Indeed, a proof of the $T(1)$ theorem can be given by using the Haar wavelets only, and these have been available since Christ's construction of his dyadic cubes with the small boundary property. In fact, Christ even formulates the general $T(1)$ theorem \cite[Theorem 8]{Christ}, but attributes it to \cite{DJS}, and proceeds to use it as a tool for proving a certain `local' variant. However, Christ's techniques would have clearly delivered a proof of the `global' $T(1)$ theorem as well, without the restrictions imposed in \cite{DJS}.

So the regular wavelets are not strictly necessary for obtaining the $T(1)$ theorem, but they nevertheless provide a rather efficient tool for that purpose. We only indicate the argument, which largely imitates the treatment in Euclidean cases given by Meyer \cite{M2}.
}

We take as space of test functions the space
 $\mathcal{V}_{s}= C^s_{0}(X)$ of functions with bounded support and  H\"older-regularity $s$ equipped with the usual topology, where $s\in (0,\eta)$ is arbitrary and we recall that $\eta$ is the regularity of the splines.  The space $\mathcal{V}_{s}$ is dense in $L^2(\mu)$ by Proposition~\ref{prop:density}.  Let $\mathcal{V}_{s}'$ denotes its dual space.  Recall the standard definition. 
 
\begin{definition}\label{def:CZK}
Let $T:\mathcal{V}_{s}\to \mathcal{V}_{s}'$ be a linear continuous operator. We say that $T$ is associated to a Calder\'on-Zygmund kernel of order $s$ if the  distributional kernel $K(x,y)$ of $T$ satisfies for some constant $C_{1}<\infty$,
$$
 |K(x,y)| \lesssim C_{1 }V(x,y)^{-1}
 $$
 when $x\ne y$ and $$
 |K(x,y)-K(x,y')|\le   C_{1}\left( \frac{d(y,y')}{d(x,y)}\right)^s  V(x,y)^{-1}
   $$
   when $x,y,y' \in X$ with $0<d(x,x') \le d(x,y)/(2A_{0})$
   $$
  |K(x,y)-K(x,y')|\le    C_{1}\left( \frac{d(x,x')}{d(x,y)}\right)^s  V(x,y)^{-1}$$
  when $x,x',y \in X$ with $0< d(x,x') \le d(x,y)/(2A_{0})$ and if furthermore, 
 for any $f \in \mathcal{V}_{s}$,  one has the representation 
  \begin{equation}
\label{eq:rep}
Tf(x)= \int K(x,y)f(y)\, d\mu(y)
\end{equation}
 for almost every $x\notin $ supp$\, f$. 
\end{definition} 

Recall also the weak boundedness property and the meaning of $T(1)$. A linear continuous operator  $T:\mathcal{V}_{s}\to \mathcal{V}_{s}'$ has weak boundedness property WBP$(\sigma)$ if  $$
 | (Tf,g)| \le C_{0} V(x,r)
 $$
 whenever $f,g \in \mathcal{V}_{\sigma}$, with support in $B(x,r)$ and are normalized by $\|f\|_{\infty}+ r^\sigma\|f\|_{\dot C^\sigma}\le 1$ and similarly for $g$. It is well-known that $WBP(\sigma)$ and $WBP(\sigma')$ are equivalent conditions whenever $s\le \sigma,\sigma'\le \eta$ when $T$ is associated to a Calder\'on-Zygmund kernel of order $s$. 
As for $T(1)$, it is defined as a continuous linear functional on the subspace of $\mathcal{V}_{s}$ of functions  $f$ with mean value 0 by 
 $$
 (T(1), f)= (Tg, f) + \int_{X}(1-g(x)) \, (^tTf)(x)\, d\mu(x)
 $$
 with $g$ a function in $\mathcal{V}_{s}$ that is 1 on a ball $B(x_{0},r)$ containing the support of $f$ and 0 on the complement of $B(x_{0},2A_{0}r)$ as in Corollary \ref{cor:existBump}, as ${}^tTf$ is integrable away  the support of $f$.

 \begin{theorem} Let $(X,d,\mu)$ be any space of homogeneous type and $T$ be associated to a Calder\'on-Zygmund kernel of order $s$. Then $T$  has a bounded extension to $L^2(\mu)$  if and only if $T$ has $WBP(s)$, $T(1) \in \BMO(\mu), {}^tT(1)\in \BMO(\mu)$. 
 \end{theorem}

 As usual it suffices to prove the converse. {With our spline wavelets at hand,} any of the standard {wavelet} proofs of the $T(1)$ theorem applies to our  statement. For example,  one can follow almost line by line the wavelet proof given in \cite[pp. 267-270]{M2}.  The weak boundedness property and the kernel representation allows to define $(T(1), \psi^k_{\alpha})$, whose absolute value  has a bound  $C\sqrt{\mu(B(y^k_{\alpha}, \delta ^k))}$.
 Then, one takes away the paraproducts  and set $U=T-\Pi_{T(1)}- \,^t\Pi_{^tT(1)}$ where
 $$
 \Pi_{\beta}f= \sum_{(k,\alpha)} (f,\bar s^{k+1}_{\alpha}) (\beta, \psi^k_{\alpha}) \psi^k_{\alpha}
 $$
 with $\bar s^{k+1}_{\alpha}$  the spline $s^{k+1}_{\alpha}$ normalized in $L^1(\mu)$. Recall that $\psi^k_{\alpha}$ is localized near the point $y^k_{\alpha}= x^{k+1}_{\alpha}\in \mathscr{Y}^{k}=\mathscr{X}^{k+1}\setminus \mathscr{X}^{k}$, 
 so that the Calder\'on-Zygmund estimates follow from the results in Section \ref{sec:technical} (see the initial comment there) when $\abs{(\beta, \psi^k_{\alpha})} \le C\sqrt{\mu(B(y^k_{\alpha}, \delta ^k))}$. The paraproduct is classically bounded on $L^2(\mu)$ if and only if the coefficients $(\beta, \psi^k_{\alpha})$ form a Carleson sequence.  It uses the maximal function $f^*(x)=\sup |(f,\bar s^{k+1}_{\alpha})|$ taken over all $k\in \Z, x^{k+1}_{\alpha}\in \mathscr{X}^{k+1}$ such that $x\in B(x^{k+1}_{\alpha}, C\delta ^{k+1})$ which is clearly comparable to Hardy-Littlewood maximal function. By Theorem \ref{th:BMO}, this is equivalent to $\beta\in \BMO(\mu)/\C$.
Thus $U$ is associated to a Calder\'on-Zygmund kernel  of order $s$ with, moreover, $U(1)=\,^tU(1)=0$. 
 
 Next, by taking $f,g$ as finite linear combination of the wavelets (which are $L^2$ dense), it suffices to prove that the coefficients of $U$ on the wavelet basis form a bounded matrix on $\ell^2$. Here, using only kernel regularity estimates (not the size, in fact), weak boundedness, $U(1)=\,^tU(1)=0$, one makes sense of the coefficients $(U\psi^k_{\alpha}, \psi^\ell_{\beta})$ and  finds the estimate
 \begin{equation}\label{eq:Ucoeff}
 |(U\psi^k_{\alpha}, \psi^\ell_{\beta})| \lesssim  \frac{C_{0}\delta ^{|k-\ell|\varepsilon}(1+\delta ^{-k\wedge \ell} d(y^{k}_{\alpha}, y^{\ell}_{\beta}))^{-\varepsilon}\sqrt{\mu(B(y^{k}_{\alpha},\delta ^{k}))\mu(B(y^{\ell}_{\beta},\delta ^{\ell}))}}{\mu(B(y^{k}_{\alpha}, \delta^{k}))+\mu(B(y^{\ell}_{\beta}, \delta ^{\ell}))+V(y^{k}_{\alpha}, y^{\ell}_{\beta})}
\end{equation}
with $0<\varepsilon<s$ and $k\wedge \ell= \inf(k,\ell)$.  The version of the Schur lemma for the $\ell^2$ boundedness
 is
 $$
 \sum_{(k,\alpha)}  |(U\psi^k_{\alpha}, \psi^\ell_{\beta})|\sqrt{\mu(B(y^{k}_{\alpha},\delta ^{k}))} \lesssim \sqrt{\mu(B(y^{\ell}_{\beta},\delta ^{\ell}))}
 $$
 uniformly in $(\ell,\beta)$ and the symmetric estimate reversing the roles of $(k,\alpha), (\ell, \beta)$. 
Details are left to the reader.

Strictly speaking, this argument works  when $\mu(X)=\infty$. When $\mu(X)<\infty$,  one has to incorporate the following observation.
 {First, the assumption that $T(1)\in\BMO(\mu)$ implies in particular that $T(1) $ is (identified to) a locally integrable function.
Second, since $X$ is bounded, it can be regarded as a ball, and the constant function $1$ is a smooth bump function associated with this ball.
Hence, the weak boundedness property implies that $\vert\int_{X} T(1)\, d\mu\vert =\vert\langle T(1),1\rangle\vert \lesssim \mu(X)$, and thus
 $\|T(1)\|_{\BMO(\mu)} +\abs{\int_{X} T(1)\, d\mu}<\infty$ (recall that our $\BMO$ norm is the homogeneous norm).} By John-Nirenberg's inequality, this implies that   $T(1) \in  L^2(\mu)$ with $\|T(1)\|_{2}\lesssim \|T(1)\|_{\BMO(\mu)} +\abs{\int_{X} T(1)\, d\mu}$.  Similarly $\|{}^tT(1)\|_{2}\lesssim  \|{}^tT(1)\|_{\BMO(\mu)} +\abs{\int_{X} {}^tT(1)\, d\mu}<\infty$. The argument before can be repeated and implies  that $\pi T\pi $ is  bounded where $\pi$ is the orthogonal projection onto  the subspace of  functions in $ L^2(\mu)$ with mean value 0 since the wavelets span this space.  The boundedness of $T$ on $L^2(\mu)$ follows readily.

\begin{remark}
As in \cite{M2}, estimate \eqref{eq:Ucoeff} is stable under matrix multiplication up to changing $\varepsilon$ to a smaller value. This algebra property furnishes a proof that  Calder\'on--Zygmund operators with $T(1)={}^tT(1)=0$ is an algebra for the composition.  This was proved in \cite{HL} by working with a new  H\"older-continuous quasi-distance as in \cite{DJS}, hence changing the class of singular integrals as discussed in the Introduction.  \end{remark}

\section{Redundancy of the size estimate}

It turns out that the size condition on the kernel is actually redundant in the $T(1)$ theorem, in that it already follows from regularity of the kernel and the weak boundedness property. To our knowledge, this remark seems new even in the context of $\R^d$ with the Lebesgue measure.
We came across this observation noticing that the regularity estimate of kernels like 
$\sum h^k_{\alpha}(x)h^k_{\alpha}(y)$ with $h^k_{\alpha}$ satisfying the size and regularity of spline or spline-wavelets is straightforward, while the size estimate required the analysis  in Section \ref{sec:technical} 
based on the structure of the point sets $\mathscr{Y}^{k}$. It is also possible to develop the $T(1)$ theory without using the size estimate at all.

To be more precise, let us say that $T:\mathcal{V}_s\to\mathcal{V}_s'$ is associated to a Calder\'on--Zygmund kernel of order $s$ \emph{in the relaxed sense}, if in Definition~\ref{def:CZK} we leave out the condition that $\abs{K(x,y)}\leq C_1 V(x,y)^{-1}$, and only assume that $K(x,y)$ is locally integrable away from the diagonal. This condition could be further relaxed, by not assuming the a priori existence of a measurable kernel at all, only that $T$ is a weak limit of operators $T_n$, which are associated to relaxed kernels of order $s$ in a uniform way. However, we stick to the stated relaxation, for the simplicity of formulation.

\begin{proposition}
Suppose that $T:\mathcal{V}_s\to\mathcal{V}_s'$ is associated to a Calder\'on--Zygmund kernel of order $s$ in the relaxed sense, and that $T$ satisfies $WBP(\sigma)$ for some $\sigma\in[s,\eta]$. Then $T$ is actually associated to a Calder\'on--Zygmund kernel of order $s$ in the sense of Definition~\ref{def:CZK}.
\end{proposition}

\begin{proof}
Recall that the paraproducts $\Pi_{T(1)}$ and ${}^t\Pi_{{}^tT(1)}$ are associated to Calder\'on-Zygmund kernels of order $s$ under these assumptions. Thus it suffices to prove the claim for the operator $U=T-\Pi_{T(1)}-{}^t\Pi_{{}^tT(1)}$ in place of $T$. The kernel of $U$ is given by
\begin{equation*}
  K(x,y)=\sum_{k,\ell,\alpha,\beta}(U\psi^k_\alpha,\psi^\ell_\beta)\psi^k_\alpha(y)\psi^\ell_\beta(x),
\end{equation*}
where the coefficients satisfy the estimate \eqref{eq:Ucoeff}, as this estimate did not make use of the size of the kernel.

We show that this series converges absolutely and satisfies the size estimate.  By symmetry, it suffices to consider the half of the sum with $\ell\leq k$, hence $\delta^\ell\geq\delta^k$, and
\begin{equation*}
\begin{split}
  &\abs{(U\psi^k_\alpha,\psi^\ell_\beta)\psi^k_\alpha(y)\psi^\ell_\beta(x)} \\
  &\qquad\lesssim\frac{\delta^{(k-\ell)\varepsilon}(1+\delta^{-\ell}d(y^k_\alpha,x^{\ell}_\beta))^{-\varepsilon}}{
     V(y^{\ell}_\beta,\delta^{\ell})+V(y^k_\alpha,y^{\ell}_\beta)}
     \exp\big(-\gamma(\delta^{-k}d(y,y^k_\alpha))^a\, \big) \exp\big(-\gamma(\delta^{-k}d(x,y^\ell_\beta))^a\, \big) \\
  &\qquad\lesssim\frac{\delta^{(k-\ell)\varepsilon}(1+\delta^{-\ell}d(y,x))^{-\varepsilon}}{
     V(x,\delta^{\ell})+V(y,x)}
     \exp\big(-\gamma(\delta^{-k}d(y,y^k_\alpha))^a\, \big) \exp\big(-\gamma(\delta^{-\ell}d(x,y^\ell_\beta))^a\, \big)
\end{split}
\end{equation*}
Hence  
\begin{equation*}
\begin{split}
  \sum_{k,\ell:\ell\leq k} &\sum_{\alpha,\beta}\abs{(U\psi^k_\alpha,\psi^\ell_\beta)\psi^k_\alpha(y)\psi^\ell_\beta(x)} \\
  &\lesssim\sum_{k,\ell:\ell\leq k}\frac{\delta^{(k-\ell)\varepsilon}(1+\delta^{-\ell}d(y,x))^{-\varepsilon}}{
     V(x,\delta^{\ell})+V(y,x)}
     \exp\big(-\gamma(\delta^{-k}d(y,\mathscr{Y}^k))^a\, \big) \exp\big(-\gamma(\delta^{-\ell}d(x,\mathscr{Y}^\ell))^a\, \big) \\
  &\lesssim\sum_{\ell:\delta^\ell\geq d(x,y)}\sum_{k:\delta^k\leq\delta^\ell}
  \frac{\delta^{(k-\ell)s}}{V(x,\delta^{\ell})} \exp\big(-\gamma(\delta^{-\ell}d(x,\mathscr{Y}^\ell))^a\, \big)
  +\sum_{k,\ell:d(x,y)\geq\delta^{\ell}\geq \delta^k}\frac{\delta^{ks}d(x,y)^{-s}}{V(x,y)},
\end{split}
\end{equation*}
where $(1+\delta^{-\ell}d(x,y))^{-s}$ was dominated by $1$ in the first sum, and by $(\delta^{-\ell}d(x,y))^{-s}=\delta^{\ell s}d(x,y)^{-s}$ in the second. In the first part, we sum a geometric series in $k$, and then use the bound of Lemma~\ref{lem:sumLargeBalls}.
In the second part, we simply sum up a geometric series in $k$ producing the factor $\delta^{\ell s}$, and then a geometric series in $\ell$.
\end{proof}

\appendix

\section{Proof of Proposition~\ref{prop:uniquenessBMO}}\label{app:BMO}

 We begin with a technical estimate:
 
\begin{lemma}\label{lem:ballComplement}
For a any ball $B=B(x_0,r)$ and $b\in\BMO(\mu)$ of norm one, we have
\begin{equation*}
  \abs{(\psi^k_\alpha,1_{B^c}(b-b_B))}\lesssim\sqrt{\mu(B^k_\alpha)}\times
  \begin{cases} 1+\log\big((\delta^{k}+r+d(y^k_\alpha,x_0))/\min\{\delta^k,r\}\big) \quad \text{always}, \\
  e^{-\gamma(r/\delta^k)^a} \qquad \text{if }y^k_\alpha\in\frac{1}{2A_0}B\text{ and }\delta^k\leq r,
  \end{cases}
\end{equation*}
with $B^k_\alpha:=B(y^k_\alpha,\delta^k)$.
\end{lemma}

\begin{proof}

By estimating $\abs{\psi^k_\alpha}\leq C\mu(B^k_\alpha)^{-1/2}e^{-\gamma\delta^{-ma}}$ on $\delta^{-m}B^k_\alpha\setminus\delta^{-(m-1)}B^k_\alpha$, and simply ignoring the indicator, we estimate up by
\begin{equation*}
 (\abs{\psi^k_\alpha},\abs{b-b_B})\leq \frac{C}{\sqrt{\mu(B^k_\alpha)}}\sum_{m=0}^{\infty}e^{-\gamma\delta^{-ma}}\int_{\delta^{-m}B^k_\alpha}\abs{b-b_B}\ud\mu,
\end{equation*}
where
\begin{equation*}
\begin{split}
  \int_{\delta^{-m}B^k_\alpha}\abs{b-b_B}\ud\mu
  &\leq \int_{\delta^{-m}B^k_\alpha}\abs{b-b_{\delta^{-m}B^k_\alpha}}\ud\mu+\mu(\delta^{-m}B^k_\alpha)\abs{b_{\delta^{-m}B^k_\alpha}-b_B} \\
  &\lesssim \mu(\delta^{-m}B^k_\alpha)\Big(1+\log\frac{\delta^{-m}\delta^k+r+d(y^k_\alpha,x_0)}{\min\{\delta^{-m}\delta^k,r)\}}\Big) \\
  &\lesssim \delta^{-mM}\mu(B^k_\alpha)\Big(1+m+\log\frac{\delta^{k}+r+d(y^k_\alpha,x_0)}{\min\{\delta^k,r\}}\Big),
\end{split}
\end{equation*}
and hence
\begin{equation*}
  (\abs{\psi^k_\alpha},\abs{b-b_B})\leq C\sqrt{\mu(B^k_\alpha)}\Big(1+\log\frac{\delta^{k}+r+d(y^k_\alpha,x_0)}{\min\{\delta^k,r\}}\Big).
\end{equation*}

If $y^k_\alpha\in\frac{1}{2A_0}B$, so that $d(y^k_\alpha,B^c)\gtrsim r$, and $\delta^k\leq r$, we make a similar estimate but exploiting the fact that for $\delta^{-m}B^k_\alpha$ to touch $B^c$, we need $\delta^{-m}\delta^k\gtrsim r(B)$. Thus
\begin{equation*}
\begin{split}
  (\abs{\psi^k_\alpha},1_{B^c}\abs{b-b_B})
   &\leq \frac{C}{\sqrt{\mu(B^k_\alpha)}}\sum_{m:\delta^{-m}\gtrsim r/\delta^k}e^{-\gamma\delta^{-ma}}\int_{\delta^{-m}B^k_\alpha}\abs{b-b_B}\ud\mu, \\
   &\leq \frac{C}{\sqrt{\mu(B^k_\alpha)}}\sum_{m:\delta^{-m}\gtrsim r/\delta^k}e^{-\gamma\delta^{-ma}}\mu(\delta^{-m}B^k_\alpha)\Big(1+m+\log\frac{r}{\delta^k}\Big) \\
   &\leq C\sqrt{\mu(B^k_\alpha)}\cdot e^{-\gamma (r/\delta^k)^a}.\qedhere
\end{split}
\end{equation*}
\end{proof}

We are ready for the proof of Proposition~\ref{prop:uniquenessBMO}, which we recall here:

\begin{proposition}
Suppose that $b\in\BMO(\mu)$ satisfies $(\psi^k_\alpha,b)=0$ for all $k,\alpha$. Then $b$ is equal to a constant.
\end{proposition}

\begin{proof}
We show for every ball $B=B(x_0,r)$ and $\varepsilon>0$ the following: there exists a constant $c$ such that $\Norm{1_B(b-c)}{\infty}\leq\varepsilon$. This clearly suffices. We may assume for simplicity that $\Norm{b}{\BMO(\mu)}\leq 1$.

Given $B$, we take a large auxiliary $\tilde{B}=B(x_0,\tilde{r})$. We use the fact that $1_{\tilde{B}}(b-b_{\tilde{B}})\in L^2(\mu)$ can be expanded in terms of the wavelets $\psi^k_\alpha$:
\begin{equation}\label{eq:1BbExpand}
\begin{split}
  1_B b
  &=1_B 1_{\tilde{B}}(b-b_{\tilde{B}})+1_B b_{\tilde{B}}
  =1_B \sum_{k,\alpha}\psi^k_\alpha(\psi^k_\alpha,1_{\tilde{B}}(b-b_{\tilde{B}}))+1_B b_{\tilde{B}} \\
  &=1_B \sum_{k:\delta^k\leq r'}\sum_{\alpha}\psi^k_\alpha(\psi^k_\alpha,1_{\tilde{B}}(b-b_{\tilde{B}}))
     +1_B \sum_{k:\delta^k> r'}\sum_{\alpha}(\psi^k_\alpha-\psi^k_\alpha(x_0))(\psi^k_\alpha,1_{\tilde{B}}(b-b_{\tilde{B}})) \\
   &\qquad  +1_B \sum_{k:\delta^k> r'}\sum_{\alpha}\psi^k_\alpha(x_0)(\psi^k_\alpha,1_{\tilde{B}}(b-b_{\tilde{B}}))
     +1_B b_{\tilde{B}}. \\
\end{split}
\end{equation}
Provided that everything converges (which we check in a moment), the terms on the last line give $1_B$ times a constant, so it suffices to show that the second-to-last line becomes arbitrarily small for properly chosen $\tilde{r}>r'>r$. For the convergence of the last line, note that
\begin{equation*}
\begin{split}
  \abs{\psi^k_\alpha(x_0)(\psi^k_\alpha,1_{\tilde{B}}(b-b_{\tilde{B}}))}
  &\leq\abs{\psi^k_\alpha(x_0)}\Norm{\psi^k_\alpha}{\infty}\Norm{1_{\tilde{B}}(b-b_{\tilde{B}})}{1} \\
  &\leq\frac{C}{\mu(B^k_\alpha)}\exp\big(-\gamma(\delta ^{-k}{d(x_0,y^k_{\alpha})})^a\, \big)
\mu(\tilde{B}),
\end{split}
\end{equation*}
and, by Lemma~\ref{lem:sumLargeBalls},
\begin{equation*}
\begin{split}
  \sum_{k:\delta^k> r'} &\sum_{\alpha}\frac{C}{\mu(B^k_\alpha)}\exp\big(-\gamma(\delta ^{-k}{d(x_0,y^k_{\alpha})})^a\, \big) \\
  &\leq\sum_{k:\delta^k>r'}\frac{C}{V(x_0,\delta^k)}\exp\big(-\gamma(\delta ^{-k}{d(x_0,\mathscr{Y}^k)})^a\, \big)
  \leq\frac{C}{V(x_0,r')}.
\end{split}
\end{equation*}
So the last line of \eqref{eq:1BbExpand} is a well-defined constant, as claimed, and it remains to see that the rest of the right side of the same equation is small.

\subsubsection*{Part $\delta^k\leq r'$.}
Since $\psi^k_\alpha$ is orthogonal to constants always, and to $b$ by assumption, we have
\begin{equation*}
  (\psi^k_\alpha,1_{\tilde{B}}(b-b_{\tilde{B}}))
  =-(\psi^k_\alpha,1_{\tilde{B}^c}(b-b_{\tilde{B}})),
\end{equation*}
and we may apply Lemma~\ref{lem:ballComplement} (with $\tilde{B}$ in place of $B$) to estimate. If $y^k_\alpha\in (\frac{1}{2A_0}\tilde{B})^c$, then
\begin{equation*}
  \abs{(\psi^k_\alpha,1_{\tilde{B}^c}(b-b_{\tilde{B}}))}
  \lesssim\sqrt{\mu(B^k_\alpha)}\Big(1+\log\frac{d(y^k_\alpha,x_0)}{\delta^k}\Big),
\end{equation*}
whereas
\begin{equation*}
  \Norm{1_B\psi^k_\alpha}{\infty} \lesssim\frac{1}{\sqrt{\mu(B^k_\alpha)}}\exp\big(-\gamma(\delta ^{-k}{d(x_0,y^k_{\alpha})})^a\, \big)
.
\end{equation*}
Hence
\begin{equation}\label{eq:smallScaleOut}
\begin{split}
  \sum_{k:\delta^k\leq r'}\sum_{\alpha:y^k_\alpha\in(\frac{1}{2A_0}\tilde{B})^c}\Norm{1_B\psi^k_\alpha(\psi^k_\alpha,1_{\tilde{B}}(b-b_{\tilde{B}}))}{\infty}
  &\lesssim  \sum_{k:\delta^k\leq r'}\sum_{\alpha:y^k_\alpha\in(\frac{1}{2A_0}\tilde{B})^c}\exp\big(-\gamma(\delta ^{-k}{d(x_0,y^k_{\alpha})})^a\, \big)
 \\
  &\lesssim  \sum_{k:\delta^k\leq r'}e^{-\gamma(\tilde{r}/\delta^k)^a}
  \lesssim e^{-\gamma(\tilde{r}/r')^a}.
\end{split}
\end{equation}
If, on the other hand, $y^k_\alpha\in\frac{1}{2A_0}\tilde{B}$, then we have the estimates
\begin{equation*}
  \abs{(\psi^k_\alpha,1_{\tilde{B}^c}(b-b_{\tilde{B}}))}
  \lesssim\sqrt{\mu(B^k_\alpha)}\cdot e^{-\gamma(\tilde{r}/\delta^k)^a},\qquad
  \Norm{1_B\psi^k_\alpha}{\infty} \lesssim\frac{1}{\sqrt{\mu(B^k_\alpha)}},
\end{equation*}
hence
\begin{equation}\label{eq:smallScaleIn}
\begin{split}
  \sum_{k:\delta^k\leq r'}\sum_{\alpha:y^k_\alpha\in\frac{1}{2A_0}\tilde{B}}\Norm{1_B\psi^k_\alpha(\psi^k_\alpha,1_{\tilde{B}}(b-b_{\tilde{B}}))}{\infty}
  &\lesssim  \sum_{k:\delta^k\leq r'}\sum_{\alpha:y^k_\alpha\in \frac{1}{2A_0}\tilde{B}} e^{-\gamma(\tilde{r}/\delta^k)^a} \\
  &\lesssim  \sum_{k:\delta^k\leq r'}e^{-\gamma(\tilde{r}/\delta^k)^a}(\tilde{r}/\delta^k)^M
  \lesssim e^{-\gamma(\tilde{r}/r')^a}.
\end{split}
\end{equation}

\subsubsection*{Part $\delta^k> r'$.}
As before, we have from Lemma~\ref{lem:ballComplement} (with $\tilde{B}$ in place of $B$) that
\begin{equation*}
\begin{split}
  \abs{(\psi^k_\alpha,1_{\tilde{B}}(b-b_{\tilde{B}}))}
  &\lesssim\sqrt{\mu(B^k_\alpha)}\Big(1+\log\frac{\delta^k+\tilde{r}+d(y^k_\alpha,x_0)}{\min\{\delta^k,\tilde{r}\}}\Big) \\
  &\lesssim\sqrt{\mu(B^k_\alpha)}\Big(1+\log_+\frac{\delta^k}{\tilde{r}}+\log_+\frac{\tilde{r}}{\delta^k}+\log_+\frac{d(y^k_\alpha,x_0)}{\delta^k}\Big),
\end{split}
\end{equation*}
and from the regularity of the wavelets that
\begin{equation*}
  \Norm{1_B(\psi^k_\alpha-\psi^k_\alpha(x_0))}{\infty}
  \lesssim\frac{1}{\sqrt{\mu(B^k_\alpha)}}\Big(\frac{r}{\delta^k}\Big)^{\eta}\exp\big(-\gamma(\delta ^{-k}{d(y^k_\alpha,x_0)})^a \, \big).
\end{equation*}
Hence
\begin{equation}\label{eq:largeScale}
\begin{split}
  \sum_{k:\delta^k> r'}&\sum_{\alpha}\Norm{1_B(\psi^k_\alpha-\psi^k_\alpha(x_0))(\psi^k_\alpha,1_{\tilde{B}}(b-b_{\tilde{B}}))}{\infty} \\
  &\lesssim  \sum_{k:\delta^k> r'}\Big(\frac{r}{\delta^k}\Big)^{\eta}\Big(1+\log_+\frac{\tilde{r}}{\delta^k}+\log_+\frac{\delta^k}{\tilde{r}}\Big) \\
  &\qquad\times\sum_{\alpha}(1+\log_+\frac{d(y^k_\alpha,x_0)}{\delta^k})\exp\big(-\gamma(\delta ^{-k}{d(y^k_\alpha,x_0)})^a \, \big).
   \\
  &\lesssim \Big(\frac{r}{r'}\Big)^{\eta}\Big(1+\log_+\frac{\tilde{r}}{r'}\Big)
\end{split}
\end{equation}

\subsubsection*{Conclusion}
Substituing \eqref{eq:smallScaleOut}, \eqref{eq:smallScaleIn} and \eqref{eq:largeScale} into \eqref{eq:1BbExpand}, we have shown that: for any $B=B(x_0,r)$, and $\tilde{r}>r'>r$, there exists a constant $c=c(B,\tilde{r},r')$ such that
\begin{equation*}
  \Norm{1_B(b-c)}{\infty}\lesssim e^{-\gamma(\tilde{r}/r')^a}+\Big(\frac{r}{r'}\Big)^{\eta}\Big(1+\log\frac{\tilde{r}}{r'}\Big).
\end{equation*}
If we now choose $r'=tr$, $\tilde{r}=tr'=t^2 r$, we get
\begin{equation*}
  \Norm{1_B(b-c_{B,t})}{\infty}\lesssim e^{-\gamma t^a}+t^{-\eta}(1+\log t)\underset{t\to\infty}{\longrightarrow} 0.
\end{equation*}
This concludes the proof.
\end{proof}
  
\bibliographystyle{acm}

\end{document}